\documentclass[11pt,a4paper]{article}
\usepackage[left=2.54cm, right=2.54cm, top=2.54cm, bottom=2.38cm]{geometry} 

\usepackage{comment}
\usepackage{amssymb}
\usepackage{amsmath}
\usepackage{amsthm} 
\usepackage{mathtools}
\usepackage[noadjust]{cite}

\usepackage[utf8]{inputenc}
\usepackage{tcolorbox}
\usepackage{ragged2e}
\usepackage{nccmath} 
\usepackage{hyperref}
\usepackage{capt-of}
\usepackage{afterpage}
\usepackage{fontenc}

\usepackage{array}    
\usepackage{booktabs} 

\usepackage{graphicx}
\usepackage{caption}
\usepackage{subcaption}	

\usepackage{accents}

\newcommand{\Hcal}{\mathcal{H}}
\newcommand{\Gcal}{\mathcal{G}}
\newcommand{\Rbb}{\mathbb{R}}
\newcommand{\Nbb}{\mathbb{N}}

\newcommand\inner[2]{\langle #1, #2 \rangle} 

\makeatletter
\newcommand{\leqnomode}{\tagsleft@true\let\veqno\@@leqno} 

\newtheorem{theorem}{Theorem}

\newtheorem{lemma}[theorem]{Lemma}
\newtheorem{remark}[theorem]{Remark}
\newtheorem{corollary}[theorem]{Corollary}

\title{The forward-backward-forward algorithm with extrapolation from the past and penalty scheme for solving monotone inclusion problems and applications}
\author{Buris Tongnoi\footnote{Faculty of Mathematics, University of Vienna, Oskar-Morgenstern-Platz 1,
Vienna 1090, Austria.; e-mail: buris.tongnoi@univie.ac.at}}
\date{}

\begin{document}
	\maketitle
	
	\textbf{Abstract:} In this paper, we propose an improved iterative method for solving the monotone inclusion problem in the form of $0 \in Ax + Dx + N_{C}(x)$ in real Hilbert space, where $A$ is a maximally monotone operator, $D$ and $B$ are monotone and Lipschitz continuous, and $C$ is the nonempty set of zeros of the operator $B$. Our investigated method, called Tseng's forward-backward-forward with extrapolation from the past and penalty scheme, extends the one proposed by Bo\textcommabelow{t} and Csetnek [Set-Valued Var. Anal. 22: 313--331, 2014]. We investigate the weak ergodic and strong  convergence (when $A$ is strongly monotone) of the iterates produced by our proposed scheme. We show that the algorithmic scheme can also be applied to minimax problems. Furthermore, we discuss how to apply the method to the inclusion problem involving a finite sum of compositions of linear continuous operators by using the product space approach and employ it for convex minimization. Finally, we present a numerical experiment in TV-based image inpainting to validate the proposed theoretical theorem.

	\textbf{keywords:} Tseng's algorithm, iterative methods, optimization problems, monotone inclusion problems, minimax problems, convergence analysis, Fitzpatrick function.
	
	\text{MSC(2010):}  	47H05, 47N10, 65K05

\section{Introdunction} \label{Motivation}

In the last decade, penalty schemes have become a popular approach for studying constrained or hierarchical optimization problems in Hilbert spaces, particularly for solving complex problems (see [\cite{Attouch2011Czarnecki,AttouchCzarnecki2011,BotCsetnek2014,2014BotCsetnek,BotCsetnekNimana2018,Noun2013Peypouquet,BanertBot2015}]). In this work, we intend to study the general monotone inclusion problem in the following form:
	\begin{align} \label{sum of three monotone operators}
		0 \in Ax + Dx + N_{C}(x),
	\end{align}
	where $A: \mathcal{H} \rightrightarrows \mathcal{H}$ is a maximally monotone operator, $D:\mathcal{H}\rightarrow\mathcal{H}$ a (single-valued) monotone and $\eta^{-1}$-Lipschitz continuous operator with $\eta>0$, and $N_C: \mathcal{H} \rightrightarrows \mathcal{H}$ is the normal cone of the closed convex set $C\subseteq \mathcal{H}$ which is the nonempty set of zeros of another operator $B : \mathcal{H} \rightarrow \mathcal{H}$, which is a (single-valued) $\mu^{-1}$-Lipschitz continuous operator with $\mu>0$. 
	
	The problem (\ref{sum of three monotone operators}) is the generalized version of the monotone inclusion problem
	
	\begin{align} \label{sum of two monotone in monotone inclusion}
		0 \in Ax + N_{C} (x),
	\end{align}
	where $C=\arg\min \Psi\neq \emptyset$, and $\Psi: \mathcal{H} \rightarrow \mathbb{R}$ 
	is a convex differentiable function with a Lipschitz gradient 
	 satisfying $\min \Psi = 0$, which introduces a penalization function with respect to the constraint $x\in C$. Implicit and explicit discretized methods to solve the problem in (\ref{sum of two monotone in monotone inclusion}) were proposed in [\cite{AttouchCzarnecki2011, Attouch2011Czarnecki}]. Several subsequent penalty type numerical algorithms in the literature are inspired by the continuous nonautonomous differential inclusion investigated in [\cite{Attouch2010Czarnecki}]. 
	
	In case $A$ is the convex subdifferential of a proper, convex, and lower semicontinuous function $\Phi: \mathcal{H} \rightarrow \bar{\mathbb{R}}$, then any solution of (\ref{sum of two monotone in monotone inclusion}) solves the convex minimization problem
		\begin{align}
			\min_{x\in\mathcal{H}} \{ \Phi(x) : x \in \arg\min \Psi\}.
		\end{align}
	Consequently, a solution to this convex minimization can be approximated by using the same iterative scheme as in [\cite{AttouchCzarnecki2011,Attouch2011Czarnecki}]. Futhermore, there are other methods for solving such optimization problem which are studied by several authors in the literature (see, for instance, [\cite{Noun2013Peypouquet,Peypouquet2012,BotCsetnekNimana2018,BanertBot2015}]).
	
	Normally, for the iterative algorithms for solving both monotone inclusion problems (\ref{sum of three monotone operators}) and (\ref{sum of two monotone in monotone inclusion}) in penalty scheme, weak ergodic convergence is proved (and strong convergence under the strong monotonicity of $A$) (see [\cite{AttouchCzarnecki2011, Attouch2011Czarnecki,BanertBot2015,2014BotCsetnek,BotCsetnek2014}]). In order to achieve the (ergodic) convergent results, the following hypotheses need to be assumed.
		\begin{itemize}
			\item For solving the problem (\ref{sum of two monotone in monotone inclusion}) in case  $C=\arg\min \Psi\neq \emptyset$ with the algorithm proposed by Attouch et al.  [\cite{AttouchCzarnecki2011,Attouch2011Czarnecki}], the Fenchel conjugate function of $\Psi$ (namely, $\Psi^{*}$) needs to fulfill some key hypotheses in this context, as shown below:
			
				$(H)\begin{cases}
					(i) &A+N_C  \;\mbox{is maximally monotone and} \; zer(A+N_C)\neq\emptyset;\\
					(ii) &\mbox{For every } p\in \operatorname{ran} (N_C),  \sum_{n\in\mathbb{N} } \lambda_n\beta_n \left[ \Psi^* \left(\frac{p}{\beta_n}\right)-\sigma_C \left(\frac{p}{\beta_n}\right)\right] < +\infty;\\
					(iii) &(\lambda_n)_{n\in\mathbb{N} }  \in \ell^2 \setminus \ell^1,
				\end{cases}$
			
			where $(\lambda_n)_{n\in\mathbb{N}}$ and $(\beta_n)_{n\in\mathbb{N}}$ are sequences of positive real numbers.  Note that the hypothesis $(ii)$ of  $(H)$ is satisfied when  $\sum_{n\in\mathbb{N}} \frac{\lambda_n}{\beta_n} < +\infty$ and $\Psi$ is bounded below by a multiple of the square of distances to $C$, as described in [\cite{Attouch2011Czarnecki}]. One such case is when $C=zer L$, where $L: \mathcal{H}\rightarrow\mathcal{H}$ is a linear continuous operator with closed range and $\Psi : \mathcal{H}\rightarrow \mathbb{R}$ is defined as $\Psi(x) = \frac{1}{2}\| Lx \|^2 $. Further situations in which condition $(ii)$ is fulfilled can be found in [\cite{AttouchCzarnecki2011}, Section 4.1].
			
			\item For solving the problem (\ref{sum of three monotone operators}) with a forward-backward type or Tseng's type (or forward-backward-forward) algorithm proposed by Bo\textcommabelow{t} et al. in [\cite{BanertBot2015,BotCsetnek2014,2014BotCsetnek}], the required hypotheses involve the Fitzpatrick function associated to the maximally monotone operator $B$ and reads as:
			
			$(H_{fitz})\begin{cases}
				(i) &A+N_C  \;\mbox{is maximally monotone and} \; zer(A+D+N_C)\neq\emptyset;\\
				(ii) &\mbox{For every } p\in \operatorname{ran} (N_C),  \sum_{n\in\mathbb{N}}  \lambda_n\beta_n \left[\sup\limits_{u\in C} \varphi_{B} \left(u,\frac{p}{\beta_n}\right)-\sigma_C \left(\frac{p}{\beta_n}\right)\right] < +\infty;\\
				(iii) &(\lambda_n)_{ n\in\mathbb{N}}  \in \ell^2 \setminus \ell^1.
			\end{cases}$
			
			\noindent where $(\lambda_n)_{n\in\mathbb{N}}$ and $(\beta_n)_{n\in\mathbb{N}}$ are sequences of positive real numbers. Note that when $B= \partial \Psi$ and $\Psi(x)=0$ for all $x\in C$, then by  (\ref{property of Fitzpatrick}) one can see that condition $(ii)$ of $(H)$ implies the second condition of $(H_{fitz})$, see [\cite{BotCsetnek2014,2014BotCsetnek}].
		\end{itemize} 
	
		Bo\textcommabelow{t} and Csetnek [\cite{BotCsetnek2014}] relaxed the cocoercivity of $B$ and $D$ to monotonicity and Lipschitzian. In this setting, a forward-backward-forward penalty type algorithm based on a method proposed by Tseng [\cite{Tseng2000}] is introduced. The convergence properties of this algorithm are studied under $(H_{fitz})$ and the condition $\lim\sup_{n \rightarrow + \infty} \left( \frac{\lambda_{n}\beta_n}{\mu} + \frac{\lambda_{n}}{\eta} \right)< 1$.
		
		In recent years, Tseng's forward-backward-forward algorithm has been modified by many researchers in various versions depending on the setting of the operators (see [\cite{Tongnoi2022,BöhmSedlmayerCsetnekBot2022,Malitsky2020Tam,Rakhlin2013Sridharan,RakhlinSridharan2013}], and references therein). One such modification is the Tseng's forward-backward-forward algorithm with extrapolation from the past, which is developed by using Popov's idea [\cite{Popov1980}] for the extragradient method with the extrapolated technique. This algorithm can store and reuse the extrapolated term in the next step of the iterative scheme, illustrated as
		$$
		\underbrace{\begin{cases} 
				y_n = J_{\gamma A}(x_n-\gamma B(x_n)),\\
				x_{n+1} =y_n+\gamma(B(x_n)-B(y_n)),
		\end{cases} }_{\text{Tseng's algorithm}}
		\quad\mbox{and} \quad
		\underbrace{\begin{cases} 
				y_n = J_{\gamma A}(x_n-\gamma B(y_{n-1})),\\
				x_{n+1} =y_n+\gamma(B(y_{n-1})-B(y_n)).
		\end{cases}}_{\text{Tseng's algorithm with extrapolation from the past}}
		$$
		where $\gamma$ satisfies some suitable condition for each method and $J_A$ is the resolvent operator of $A$. According to this scheme, it has the potential to reduce computational costs and energy consumption in practical applications.
		
		Motivated by these considerations, we investigate the Tseng's forward-backward-forward algorithm with extrapolation from the past involving a penalty scheme under certain hypotheses for solving the inclusion problem in (\ref{sum of three monotone operators}). We prove its weak ergodic convergence and further the strong convergence when the operator $A$ is a strongly monotone operator in Section \ref{forward-backward-forward Penalty Schemes}. Furthermore, we can extend our results in Section \ref{forward-backward-forward Penalty Schemes} to solve minimax problems as elaborated in Section \ref{Minimax optimization problem}.  In Section \ref{The FBF algorithm with extrapolation penalty scheme for the problem involving linearly composed and parallel-sum type monotone operators}, our proposed algorithm can be extended to solving more intricate problems involving the finite sum of composition of linear continuous operator with maximally monotone operators by using the product space approach, and the convergence results for our modified iterative method are also provided.
		
		To broaden the applicability of our scheme, we also provide the iterative scheme and its convergence result for the convex minimization problem expressed in Section \ref{Convex minimization problem}. Finally, we demonstrate a numerical experiment in TV-based image inpainting to ensure theoretical convergence results in Section \ref{A Numerical Experiment in TV-Based Image Inpainting}.

\section{Preliminaries}\label{Preliminaries}

In this paper, we introduce some notations used throughout. We denote the set of \textit{positive integers} by $\mathbb{N}$ and a real Hilbert space with an \textit{inner product} $\inner{\cdot}{\cdot}$ and the associated \textit{norm} $\| \cdot \| = \sqrt{\inner{\cdot}{\cdot}}$ by $\mathcal{H}$. A sequence $(x_n)_{n\in\mathbb{N}}$ is said to \textit{converge weakly} to $x$, if for any $y\in\mathcal{H}$, $\inner{x_n}{y}\rightarrow\inner{x}{y}$, and we use the symbols $\rightharpoonup$ and $\rightarrow$ to represent weak and strong convergence, respectively. For a linear continuous operator $L: \mathcal{H}\rightarrow \mathcal{G}$, where $\mathcal{G}$ is another real Hilbert space, we define the norm of $L$ as $\| L \| = \sup \{ \|Lx\| : x\in \mathcal{H}, \|x\| \leq 1 \}$. We also use $L^* : \mathcal{G} \rightarrow \mathcal{H}$ to denote the adjoint operator of $L$, defined by $\inner{L^* y}{x} = \inner{y}{Lx}$ for all $(x,y) \in \mathcal{H} \times \mathcal{G}$.

For a function $f: \mathcal{H} \rightarrow \bar{\mathbb{R}}$ (where $ \bar{\mathbb{R}}:= \mathbb{R}\cup\{\pm\infty \}$ is the extended real line), we denote its \textit{effective domain} by $\operatorname{dom} f = \{ x \in \mathcal{H} : f(x) < + \infty \}$ and say that $f$ is \textit{proper} if $\operatorname{dom} f\neq\emptyset$ and $f(x) \neq - \infty$ for all $x\in \mathcal{H}$. Let $f^* : \mathcal{H} \rightarrow \bar{\mathbb{R}}$, where $f^*(u) = \sup_{x\in\mathcal{H}} \{ \inner{u}{x} - f(x) \}$ for all $u\in\mathcal{H}$, be the \textit{conjugate function} of $f$. We denote by $\Gamma (\mathcal{H})$ the family of proper, convex, and lower semicontinuous extended real-valued functions defined on $\mathcal{H}$. The \textit{subdifferential} of $f$ at $x\in\mathcal{H}$, with $f(x)\in\mathbb{R}$, is the set $\partial f(x) := \{v \in \mathcal{H} : f(y) \geq f(x) + \inner{v}{y-x} \; \forall y\in\mathcal{H} \}$. We take, by convention, $\partial f(x) :=\emptyset$ if $f(x) \in { \pm \infty }$.

The \textit{indicator function} of a nonempty set $S\subseteq \mathcal{H}$, denoted by $\delta_{S} : \mathcal{H}\rightarrow \bar{\mathbb{R}}$, is defined as follows: $\delta_{S}(x)=0$ if $x\in S$ and $\delta_{S}(x)=+\infty$ otherwise. The subdifferential of the indicator function is called the \textit{normal cone} of $S$, denoted by $N_{S}(x)$. For $x \in S$, $N_{S}(x)=\{ u\in\mathcal{H} : \inner{y-x}{u} \leq 0 \; \forall y \in S \}$, and $N_{S}(x) = \emptyset$ if $x\notin S$. It is worth noting that for $x\in S$, $u \in N_{S} (x)$ if and only if $\sigma_{S}(u) = \inner{x}{u}$, where $\sigma_{S}$ is the support function of $S$. The \textit{support function} is defined as $\sigma_{S} (u) = \sup_{y\in S} \inner{y}{u}$. A \textit{cone} is a set that is closed under positive scaling, i.e., $\lambda x \in S$ for all $x \in S$ and $\lambda > 0$. The normal cone, which is a tool used in optimization, is always a cone itself.


Let $A$ be a set-valued operator mapping from a Hilbert space $\mathcal{H}$ to itself. We define the \textit{graph} of $A$ as $GrA=\{(x,u)\in\mathcal{H}\times\mathcal{H}:u\in Ax\}$, the \textit{domain} of $A$ as $\operatorname{dom} A=\{x\in\mathcal{H}:Ax\neq\emptyset\}$, and the \textit{range} of $A$ as $\operatorname{ran} A=\{u\in\mathcal{H}:\exists x\in\mathcal{H} \;\mbox{s.t.} \; u\in Ax\}$. The \textit{inverse operator} of $A$, denoted by $A^{-1}:\mathcal{H}\rightrightarrows \mathcal{H}$, is defined as $(u,x)\in Gr A^{-1}$ if and only if $(x,u)\in Gr A$.
We also define the \textit{set of zeros} of $A$ as $zerA=\{x\in\mathcal{H}:0\in Ax\}$. The operator $A$ is said to be \textit{monotone} if $\inner{x-y}{u-v}\geq 0$ for all $(x,u),(y,v)\in Gr A$. Moreover, a monotone operator $A$ is said to be \textit{maximally monotone} if its graph on $\mathcal{H}\times\mathcal{H}$ has no proper monotone extension. If $A$ is maximally monotone, then $zerA$ is a closed and convex set, and a point $z\in\mathcal{H}$ belongs to $zerA$ if and only if $\inner{u-z}{w}\geq 0$ for all $(u,w)\in Gr A$. Conditions ensuring that $zerA$ is nonempty can be found in [\cite{BC-Book}, Section 23.4]. 

The set-valued operator $A$ is said to be $\gamma$-\textit{strongly monotone} with $\gamma > 0$ if and only if $\inner{x-y}{u-v}\geq \gamma \| x-y\|^2$ for all $(x,u),(y,v) \in Gr A$. If $A$ is maximally monotone and strongly monotone, then the set of zeros $zer A$ is a singleton and thus nonempty (refer to [\cite{BC-Book}, Corollary 23.37]). A single-valued operator $A : \mathcal{H}\rightarrow \mathcal{H}$ is considered to be $\gamma$-\textit{cocoercive} if $\inner{x-y}{Ax-Ay} \geq \gamma \|Ax-Ay\|^2$, and $\gamma$-\textit{Lipschitz continuous} (or $\gamma$-Lipschitzian) if $\|Ax-Ay\| \leq \gamma \|x-y\|$ for all $(x,y)\in \mathcal{H}\times \mathcal{H}$.

In optimization, the proximal operator is a common tool used to develope algorithms for solving problems that involve a sum of a smooth function and a nonsmooth function. The \textit{proximal operator} of a function $f$ at a point $x$ is defined as the unique minimizer of the function
$y\mapsto f(y) + \frac{1}{2}\|y-x\|^2$, denoted by $\text{prox}_{ f}(x):\mathcal{H}\rightarrow\mathcal{H}$. 
The \textit{resolvent} of $A$, denoted by $J_A : \mathcal{H} \rightrightarrows \mathcal{H}$, is defined as $J_{A}=(Id+A)^{-1}$ where $Id: \mathcal{H} \rightarrow \mathcal{H}$, $Id(x) = x$ for all $x\in\mathcal{H}$, is the \textit{identity operator} on $\mathcal{H}$. Moreover, if $A$ is maximally monotone, then $J_A: \mathcal{H}\rightarrow \mathcal{H}$ is single-valued and maximally monotone (see [\cite{BC-Book}, Corollary 23.10]). Notice that $J_{ \partial f} = (Id +  \partial f)^{-1}=\text{prox}_{f}$ (see also [\cite{BC-Book}, Proposition 16.34]).

Fitzpatrick introduced the notation we use in this paper in [\cite{Fitzpatrick1988}], and it has proved to be a useful tool for investigating the maximality of monotone operators using convex analysis techniques (see [\cite{Bauschke2006, BC-Book}]). For a monotone operator $A$, we define its associated \textit{Fitzpatrick function} $\varphi_{A}: \mathcal{H}\times\mathcal{H}\rightarrow\bar{\mathbb{R}}$ as follows:
\begin{align*}
	\varphi_{A}(x,u) = \sup_{(y,v)\in Gr A} \{ \langle x,v\rangle + \langle y,u\rangle - \langle y,v\rangle \}.
\end{align*}
This function is convex and lower semicontinuous, and it will be an important tool for investigating convergence in this paper.

If $A$ is a maximally monotone operator, then its Fitzpatrick function is proper and satisfies $\varphi_{A}(x,u) \geq \langle x,u\rangle$ for all $(x,u) \in \mathcal{H}\times\mathcal{H}$, with equality if and only if $(x,u)\in Gr A$. In particular, the subdifferential $\partial f$ of any convex function $f\in \Gamma(H)$ is a maximally monotone operator (cf. [\cite{Rockafellar1970}]), and we have $(\partial f)^{-1} = \partial f^*$. Furthermore, we have the inequality (see [\cite{Bauschke2006}])
\begin{align}\label{property of Fitzpatrick}
	\varphi_{\partial f}(x,u) \leq f(x) + f^*(u) \quad \forall (x,u)\in \mathcal{H}\times \mathcal{H},
\end{align}
where $f^*$ is the convex conjugate of $f$. 

Before we enter on a detailed analysis of the convergence result, we introduce some definitions of sequences and useful lemmas that will be used several times in the paper. Let $(x_n)_{n\in\mathbb{N}\cup\{0\}}$ be a sequence in $\mathcal{H}$ and $(\lambda_{k})_{k\in\mathbb{N}\cup\{0\}}$ a sequence of positive numbers such that $\sum_{k\in\mathbb{N}\cup\{0\}} \lambda_k = +\infty$. Let $(z_n)_{n\in\mathbb{N}\cup\{0\}}$ be the sequence of weighted averages defined as shown in [\cite{Attouch2011Czarnecki,2014BotCsetnek,BotCsetnek2014}]:
\begin{align}\label{weigthed_averages}
	z_n = \frac{1}{\tau_n} \sum_{k=0}^{n} \lambda_k x_k,\quad \mbox{where}\; \tau_n = \sum_{k=0}^{n} \lambda_k \; \forall n \in \mathbb{N}.
\end{align}

The following lemma is the Opial-Passty Lemma, which serves as a key tool for achieving the convergence results in our analysis.

\begin{lemma}[Opial-Passty]\label{Lemma1}
	Let $T$ be a nonempty subset of real Hilbert space $\mathcal{H}$ and assume that  $\lim_{n\rightarrow \infty} \|x_n -x \|$ exists for every $x\in T$. If every weak cluster point of $(x_n)_{n\in\mathbb{N}}$(respectively $(z_n)_{n\in\mathbb{N}}$) belongs to $T$, then $(x_n)_{n\in\mathbb{N}}$(respectively $(z_n)_{n\in\mathbb{N}}$) converges weakly to an element in $T$ as $n\rightarrow \infty$.
\end{lemma}

Another key lemma underlying our work is presented below, which establishes the existence of convergence and the summability of sequences satisfying the conditions of the lemma. This lemma has been sourced from [\cite{AttouchCzarnecki2011}].

\begin{lemma}\label{Lemma2}
	Let $(a_n)_{n\in\mathbb{N}},(b_n)_{n\in\mathbb{N}}$ and $(\epsilon_{n})_{n\in\mathbb{N}}$ be real sequences. Assume that $(a_n)_{n\in\mathbb{N}}$ is bounded from below, $(b_n)_{n\in\mathbb{N}}$ is nonnegative sequences, $(\epsilon_n)_{n\in\mathbb{N}} \in \ell^{1}$ such that $$a_{n+1} \leq a_n - b_n + \epsilon_n  \quad \forall n \in \mathbb{N}$$ 
	Then $(a_n)$ is convergent and $(b_n)_{n\in\mathbb{N}} \in \ell^{1}$.
\end{lemma}

\section{Forward-backward-forward with extrapolation and penalty schemes}\label{forward-backward-forward Penalty Schemes}

In this section, we will begin by presenting a problem introduced in [\cite{BotCsetnek2014,2014BotCsetnek}], which can be formulated as follows:\\

	\noindent \textbf{Problem 1:} Let $\mathcal{H}$ be a real Hilbert space, $A : \mathcal{H}\rightrightarrows \mathcal{H}$ a maximally monotone operator, $D : \mathcal{H}\rightarrow \mathcal{H}$ a monotone and $\eta^{-1}$-Lipschitz continuous operator with $\eta >0$, $B:\mathcal{H}\rightarrow \mathcal{H}$ a monotone and $\mu^{-1}$-Lipschitz continuous operator with $\mu >0$ and suppose that $C = zerB\neq\emptyset$. The monotone inclusion problem that we want to solve is 
	\begin{align} \label{monotone inclusion problem_1}
		0\in Ax+Dx+N_{C}(x).
	\end{align}


Some methods to solve \textbf{Problem 1} have been already studied in [\cite{BotCsetnek2014,2014BotCsetnek}] by using Tseng's algorithm and penalty scheme. The iterative scheme presented below for solving \textbf{Problem 1} draws inspiration from Tseng's forward-backward-forward algorithm with extrapolation from the past (see [\cite{Tongnoi2022,BöhmSedlmayerCsetnekBot2022,Malitsky2020Tam}]). 
	
	\begin{tcolorbox}[colback=white]
		\textbf{Algorithm 1}:
		\begin{equation}
		\begin{aligned} 
			\mbox{Initialization}&:\; \mbox{Choose}\; y_{-1},x_0\in\mathcal{H} \;\mbox{with}\; y_{-1}=x_0. \\
			\mbox{For}\; n\in\mathbb{N}\cup\{ 0\}\; \mbox{set}&:\; y_n = J_{\lambda_n A}\left[ x_n-\lambda_n D(y_{n-1})-\lambda_n\beta_n B(y_{n-1})\right];\\
			&\quad x_{n+1}=\lambda_n \beta_n \left[(B(y_{n-1})-B(y_n))\right]+\lambda_n\left[D(y_{n-1})-D(y_n)\right]+y_n,
		\end{aligned}
		\end{equation}
	\end{tcolorbox}
	\noindent where $(\lambda_n)_{n\in\mathbb{N}\cup\{ 0\}}$ and $(\beta_n)_{n\in\mathbb{N}\cup\{ 0\}}$ are sequences of positive real numbers.
	
	Of course, in order to demonstrate the (ergodic) convergence results we need to assume the hypotheses $(H_{fitz})$ for every $n\in\mathbb{N}\cup\{ 0\}.$
	
	\begin{remark}\ 
		\begin{enumerate}
		\item When the operator $B$ satisfies $Bx=0$ for all $x\in\mathcal{H}$ (which implies $N_C(x)={0}$ for every $x\in\mathcal{H}$), the algorithm outlined in \textbf{Algorithm 1} reduces to the error-free scenario with the identity variable matrices of the forward-backward-forward with extrapolation method, as introduced in [\cite{Tongnoi2022}].
		\item Actually, we can choose $y_{-1}$ as any starting point in $\mathcal{H}$. 
		\end{enumerate}
	\end{remark}

Prior to establishing the weak ergodic convergence result, we shall prove the following lemma, which serves as a valuable tool for our main outcome.

\begin{lemma}\label{Lemma decreasing ineqality}
	Let $(x_n)_{n\in\mathbb{N}\cup\{0\} }$ and $(y_n)_{ { n\in\mathbb{N}\cup\{0,-1\} } }$ be the sequence generated by \textbf{Algorithm 1} and let be $(u,w)\in Gr(A+D+N_C)$ such that $w=v+p+D(u)$, where $v\in A(u)$ and $p\in N_C(u)$. Then the following inequlity holds for $n\in\mathbb{N}\cup\{0\}$:
		\begin{equation} \label{Seq.lemma}
		\begin{aligned} 
		&\|x_{n+1}-u\|^2 -\|x_n-u\|^2 +\left(\frac{1}{2}-\left(\frac{\lambda_n\beta_n}{\mu}+\frac{\lambda_n}{\eta}\right)^2\right) \|y_{n-1}-y_n\|^2\\
		&\leq \left(\frac{\lambda_{n-1}\beta_{n-1}}{\mu}+\frac{\lambda_{n-1}}{\eta}\right)^2 \|y_{n-2}-y_{n-1}\|^2 + 2 \lambda_n\beta_n \left[\sup_{u\in C}\varphi_B (u,\frac{p}{\beta_n})-\sigma_C\left(\frac{p}{\beta_n}\right) \right] + 2\lambda_n\inner{u-y_n}{w}.
		\end{aligned}
		\end{equation}
\end{lemma}

\begin{proof}
	From \textbf{Algorithm 1} and the definition of the resolvent, we can derive that $\forall n\in\mathbb{N}$
	\begin{equation}
		\begin{aligned}\label{1_equality of lemma proof}
		y_n = (Id+\lambda_n A)^{-1} [x_n-\lambda_n D(y_{n-1})-\lambda_n\beta_n B(y_{n-1})] 
		\Leftrightarrow \frac{x_n-y_n}{\lambda_n} - D(y_{n-1}) -\beta_n B(y_{n-1})\in A(y_n).
		\end{aligned}
	\end{equation}
	Thus, we have $\frac{x_n-y_n}{\lambda_n}  -\beta_n B(y_{n-1})- D(y_{n-1})\in A(y_n)$ for ever $n\in \mathbb{N}$.
	Since $v\in Au$, then we can guarantee by using the monotonicity of $A$ that 
	\begin{small}
	\begin{equation}
		\begin{aligned}
			\inner{y_n-u}{\frac{x_n-y_n}{\lambda_n}  -\beta_n B(y_{n-1})- D(y_{n-1})-v}\geq 0 \Leftrightarrow \inner{y_n-u}{x_n-y_n  -\lambda_n\beta_n B(y_{n-1})- \lambda_nD(y_{n-1})-\lambda_n v}\geq 0.
		\end{aligned}
	\end{equation}
	\end{small}
	By using the definition of $x_{n+1}$, we obtain that $\forall n\in\Nbb$,
	\begin{equation} \label{inner of u-y_n and x_n-y_n}
	\begin{aligned}
		&\inner{u-y_n}{x_n-y_n} \\
		&\leq\inner{u-y_n}{\lambda_n\beta_n B(y_{n-1}) +\lambda_n D(y_{n-1}) + \lambda_nv}\\
		&\leq \inner{u-y_n}{[\lambda_n\beta_n(B(y_{n-1})-B(y_n))+\lambda_n (D(y_{n-1})-D(y_{n}))+y_n ] + \lambda_n\beta_nB(y_{n}) + \lambda_n D(y_n)+\lambda_nv-y_n} \\
		&= \inner{u-y_n}{x_{n+1}-y_n+\lambda_n\beta_nB(y_n)+\lambda_nD(y_n)+\lambda_n v }\\
		&= \inner{u-y_n}{x_{n+1}-y_n} + \inner{u-y_n}{\lambda_n\beta_nB(y_n)} +\inner{u-y_n}{\lambda_nD(y_n)}+\inner{u-y_n}{\lambda_nv}.
	\end{aligned}
	\end{equation}
	For any $a_1,a_2,a_3 \in \Hcal$, we know that $\frac{1}{2}\|a_1-a_2\|^2 - \frac{1}{2}\|a_3-a_1 \|^2 +\frac{1}{2}\|a_3-a_2\|^2 =  \inner{a_1-a_2}{a_3-a_2}$. Then, it follows from (\ref{inner of u-y_n and x_n-y_n}), we have that 
	
		\begin{align*}
		\frac{1}{2}\|u-y_n\|^2 - \frac{1}{2}\|x_n-u \|^2 +\frac{1}{2}\|x_n-y_n\|^2 &\leq \frac{1}{2}\|u-y_n\|^2 - \frac{1}{2}\|x_{n+1}-u \|^2 +\frac{1}{2}\|x_{n+1}-y_n\|^2\\
		&\;+\inner{u-y_n}{\lambda_n\beta_nB(y_n)} +\inner{u-y_n}{\lambda_nD(y_n)}+\inner{u-y_n}{\lambda_nv}. 
		\end{align*} 
	Because $v=w-p-Du$ and the fact that $D$ is monotone, we have that for every $n\in\Nbb$
	\begin{equation} \label{ineq_sup_sigma}
	\begin{aligned}
		\|x_{n+1}-u\|^2 -\|x_{n}-u\|^2 
		&\leq \|x_{n+1}-y_n\|^2 -\|x_{n}-y_n\|^2 + 2 \lambda_n\beta_n \inner{u-y_n}{B(y_n)}+ 2\lambda_n\inner{u-y_n}{D(y_n)}\\ &\quad+2\lambda_n\inner{u-y_n}{w-p-Du}\\
		&= \|x_{n+1}-y_n\|^2 -\|x_{n}-y_n\|^2 + 2 \lambda_n\beta_n \inner{u-y_n}{B(y_n)} \\
		&\quad + 2\lambda_n\beta_n\inner{u-y_n}{\frac{-p}{\beta_n}} + 2\lambda_n\inner{u-y_n}{D(y_n)-Du} +2\lambda_n\inner{u-y_n}{w}\\
		&=  \|x_{n+1}-y_n\|^2 -\|x_{n}-y_n\|^2 \\
		&\quad + 2 \lambda_n\beta_n \!\left(\!\!\inner{u}{B(y_n)}\!+ \inner{y_n}{\frac{p}{\beta_n}}\!-\inner{y_n}{B(y_n)}\!-\inner{u}{\frac{p}{\beta_n}}\!\right)\\
		&\quad + 2\lambda_n\underbrace{\inner{u-y_n}{D(y_n)-Du}}_{\leq 0} +2\lambda_n\inner{u-y_n}{w}\\
		&\leq  \|x_{n+1}-y_n\|^2 -\|x_{n}-y_n\|^2 \\
		&\quad + 2 \lambda_n\beta_n \left[ \sup_{u\in C} \left\lbrace \sup_{(s, Bs)\in GrB}  \{\inner{u}{B(s)}\!+ \inner{s}{\frac{p}{\beta_n}}\!-\inner{s}{B(s)} \} \right\rbrace - \sigma_C\left(\frac{p}{\beta_n}\right)  \right]\\
		&\quad +2\lambda_n\inner{u-y_n}{w}\\
		&=  \|x_{n+1}-y_n\|^2 -\|x_{n}-y_n\|^2 + 2 \lambda_n\beta_n \left[ \sup_{u\in C}\varphi_{B} \left(u,\frac{p}{\beta_n}\right) - \sigma_C\left(\frac{p}{\beta_n}\right) \right] \\
		&\quad +2\lambda_n\inner{u-y_n}{w}.
	\end{aligned}
	\end{equation}
	
	\noindent By the Lipschitz continuity of $B$ and $D$, it yields
	\begin{equation}\label{diff_x_n+1_&_y_n}
	\begin{aligned}
		\|x_{n+1}-y_n\|
		&=\| \lambda_n\beta_n(B(y_{n-1})-B(y_n))+\lambda_n(D(y_{n-1})-D(y_n)) \|\\ 
		&\leq\lambda_n\beta_n\|B(y_{n-1})-B(y_n)\|+\lambda_n\| D(y_{n-1})-D(y_n) \|\\
		&= \left(\frac{\lambda_n\beta_n}{
			\mu}+\frac{\lambda_n}{
			\eta} \right) \| y_{n-1}-y_n\|.
	\end{aligned}
	\end{equation}
	It follows from (\ref{ineq_sup_sigma}) and above inequality that
	\begin{align}\label{first inequality form}
	\|x_{n+1}-u\|^2-\|x_n-u\|^2
	&\leq \left(\frac{\lambda_n\beta_n}{\mu}+\frac{\lambda_n}{\eta}\right)^2 \| y_{n-1}-y_n\|^2-\|x_n-y_n\|^2\\
	&\quad +2\lambda_n\beta_n\left[ \sup_{u\in C}\varphi_{B} \left(u,\frac{p}{\beta_n}\right) - \sigma_C\left(\frac{p}{\beta_n} \right) \right] +2\lambda_n\inner{u-y_n}{w}.
	\end{align}
	By Parallelogram law, we know that $2\|x_n-y_n\|^2+2\|x_n-y_{n-1}\|^2 = \|y_{n-1}-y_n\|^2+\|(x_n-y_n)+(x_n-y_{n-1})\|^2$. Then we have
		$$ -\|x_n-y_{n}\|^2 \leq \|x_n-y_{n-1}\|^2 -\frac{1}{2} \|y_{n-1}-y_n\|^2.$$
	So, we can derive from (\ref{first inequality form}) that
	\begin{align*}
		\|x_{n+1}-u\|^2-\|x_n-u\|^2
		&\leq \left(\frac{\lambda_n\beta_n}{\mu}+\frac{\lambda_n}{\eta}\right)^2 \| y_{n-1}-y_n\|^2+\left[ \|x_n-y_{n-1}\|^2 -\frac{1}{2} \|y_{n-1}-y_n\|^2\right]\\
		&\quad +2\lambda_n\beta_n\left[ \sup_{u\in C}\varphi_{B} \left(u,\frac{p}{\beta_n}\right) - \sigma_C\left(\frac{p}{\beta_n}  \right)\right]+2\lambda_n\inner{u-y_n}{w}.
	\end{align*}
	It follows from (\ref{diff_x_n+1_&_y_n}) that $\|x_n-y_{n-1}\| \leq \left(\frac{\lambda_{n-1}\beta_{n-1}}{\mu}+\frac{\lambda_{n-1}}{\eta}\right) \| y_{n-2}-y_{n-1}\|$.
	Thus,
	\begin{align*}
		\|x_{n+1}-u\|^2-\|x_n-u\|^2
		&\leq \left(\frac{\lambda_n\beta_n}{\mu}+\frac{\lambda_n}{\eta}\right)^2 \| y_{n-1}-y_{n}\|^2\\
		&\quad +\left[ \left(\frac{\lambda_{n-1}\beta_{n-1}}{\mu}+\frac{\lambda_{n-1}}{\eta}\right)^2 \| y_{n-2}-y_{n-1}\|^2 -\frac{1}{2} \|y_{n-1}-y_n\|^2\right]\\
		&\quad +2\lambda_n\beta_n\left[ \sup_{u\in C}\varphi_{B} \left(u,\frac{p}{\beta_n}\right) - \sigma_C\left(\frac{p}{\beta_n}  \right)\right]+2\lambda_n\inner{u-y_n}{w}.
	\end{align*}
	Therefore
	\begin{align*} 
		&\|x_{n+1}-u\|^2-\|x_n-u\|^2+\left(\frac{1}{2}-\left(\frac{\lambda_n\beta_n}{\mu}+\frac{\lambda_n}{\eta}\right)^2\right) \|y_{n-1}-y_n\|^2\\
		&\leq \left(\frac{\lambda_{n-1}\beta_{n-1}}{\mu}+\frac{\lambda_{n-1}}{\eta}\right)^2 \|y_{n-2}-y_{n-1}\|^2 + 2 \lambda_n\beta_n \left[\sup_{u\in C}\varphi_B \left(u,\frac{p}{\beta_n}\right)-\sigma_C\left(\frac{p}{\beta_n}\right) \right] + 2\lambda_n\inner{u-y_n}{w}.
	\end{align*}
\end{proof}

Subsequently, we state the convergence of \textbf{Algorithm 1} below.

\begin{theorem} \label{Thm_weak (ergodic) convergence}
	Let $(x_n)_{n\in\mathbb{N}\cup\{0\}}$ and $(y_n)_{n\in\Nbb\cup\{0,-1\}}$ be the sequences generated by \textbf{Algorithm 1} and $(z_n)_{n\in\mathbb{N}\cup\{0\}}$ be the sequence defined in (\ref{weigthed_averages}). If ($H_{fitz}$) is fulfilled and $\limsup_{n\rightarrow+\infty} \left(\frac{\lambda_n\beta_n}{
	\mu}+\frac{\lambda_n}{
	\eta} \right) < \frac{1}{2}$,
	then $(z_n)_{n\in\mathbb{N}\cup\{0\}}$ converges weakly to an element in $zer(A+D+N_C)$ as $n\rightarrow+\infty$.
\end{theorem}

\begin{proof}
	The proof of the theorem divides into three parts as the following statements:
	\begin{enumerate}
		\item[(a)] $\sum_{n\in\mathbb{N}\cup\{0\}} \| y_{n} - y_{n-1} \|^2 < +\infty$ and the sequence $(\|x_n-u\|)_{n\in\Nbb\cup\{0\}}$ is convergent for every $u\in zer(A+D+N_C)$;
		\item[(b)] every weak cluster point of $(z'_n)_{n\in\mathbb{N}\cup\{ 0 \}}$ lies in  $zer(A+D+N_C)$, where
		\begin{align*}
			z'_n=\frac{1}{\tau_n} \sum_{k=0}^{n} \lambda_ky_k \quad \mbox{and}\quad \tau_n = \sum_{k=0}^{n} \lambda_k \; \forall n\in\Nbb; 
		\end{align*}
	\item[(c)] every weak cluster point of $(z_n)_{n\in\Nbb\cup\{0\}}$ lies in  $zer(A+D+N_C)$.
	\end{enumerate}
	Every step of proof can be showed as follows.
	
	\noindent(a) Let $u\in zer(A+D+N_C)$. For $n\in\mathbb{N}$, we take $w=0$ in (\ref{Seq.lemma}), so we have
	\begin{equation}\label{1_ineq in case w=0}
	\begin{aligned}
		&\|x_{n+1}-u\|^2-\|x_n-u\|^2+\left(\frac{1}{2}-M_n^2\right) \|y_{n-1}-y_n\|^2\\
		&\leq M_{n-1}^2 \|y_{n-2}-y_{n-1}\|^2 + 2 \lambda_n\beta_n \left[\sup_{u\in C}\varphi_B \left(u,\frac{p}{\beta_n}\right)-\sigma_C\left(\frac{p}{\beta_n}\right) \right],
	\end{aligned}
	\end{equation}
	where $M_n=\frac{\lambda_n\beta_n}{\mu}+\frac{\lambda_n}{\eta}$. Rearranging the inequality,  we have that
	\begin{equation} \label{2_ineq in case w=0}
	\begin{aligned}
		&\|x_{n+1}-u\|^2 + M_n^2 \|y_{n-1}-y_n\|^2+\left(\frac{1}{2}-2M_n^2\right) \|y_{n-1}-y_n\|^2\\
		&\leq \|x_n-u\|^2 + M_{n-1}^2 \|y_{n-2}-y_{n-1}\|^2 + 2 \lambda_n\beta_n \left[\sup_{u\in C}\varphi_B \left(u,\frac{p}{\beta_n}\right)-\sigma_C\left(\frac{p}{\beta_n}\right) \right].
	\end{aligned}	
	\end{equation}
	By the hypothesis, we have that $\liminf_{n\rightarrow +\infty} \left( \frac{1}{2}- 2M_n^2 \right) >0$, it follows from Lemma \ref{Lemma2} that the sequence $E_n:= \|x_n-u\|^2 + M_{n-1}^2 \|y_{n-2}-y_{n-1}\|^2 $ converges and $\sum_{n\in\mathbb{N}}  \|y_{n-1}-y_n\|^2 < +\infty$.  This means that $(\|y_{n-1}-y_n\|^2)_{n\in\Nbb\cup\{0\}} \rightarrow 0$ and so $(\|x_n-u\|)_{n\in\Nbb\cup\{0\}}$ is a convergent sequence.
	\\

	\noindent (b)  Let $z$ be a weak cluster point of $(z'_n)_{n\in\mathbb{N}}$.
	Take $(u,w)\in Gr(A+D+N_C)$ such that $w=v+p+Du$, 	where $v\in Au$ and $p\in N_C(u)$. Let $K\in\Nbb$ 	with  $\forall n_0\in\Nbb$, $K\geq n_0+2$. Summing up for $n=n_0+1,	\dots,K$ the inequalities in (\ref{Seq.lemma})  (it is nothing else than (\ref{1_ineq in case w=0}) with the additional term $2\lambda_n\inner{u-y_n}{w}$) 
	 with $\limsup_{n\rightarrow+\infty} \left(\frac{\lambda_n\beta_n}{
		\mu}+\frac{\lambda_n}{
		\mu} \right) < \frac{1}{2}$, for $n\in\mathbb{N}$ and $M_n=\frac{\lambda_n\beta_n}{\mu}+\frac{\lambda_n}{\eta}$ we get that 
	\begin{align*}
		\sum_{n=n_0+1}^{K} \|x_{n+1}-u\|^2 - \sum_{n=n_0+1}^{K}\|x_n-u\|^2 
		&\leq \sum_{n=n_0+1}^{K} M_{n-1}^2 \| y_{n-2}-y_{n-1}\|^2 \\
		&\quad + 2 \sum_{n=n_0+1}^{K} \lambda_n\beta_n \left[\sup_{u\in C}\varphi_B \left(u,\frac{p}{\beta_n}\right)-\sigma_C\left(\frac{p}{\beta_n}\right) \right] \\
		&\quad + 2 \left\langle \sum_{n=0}^{K} \lambda_n u - \sum_{n=0}^{K} \lambda_n y_n - \sum_{n=0}^{n_0}\lambda_n u + \sum_{n=0}^{n_0} \lambda_n y_n, w \right\rangle.
	\end{align*}
	This leads to the following inequality.
	\begin{small}
		\begin{align*}	
			\|x_{K+1}-u\|^2-\|x_{n_0+1}-u\|^2 
			&\leq  \sum_{n=n_0+1}^{+\infty} M_{n-1}^2 \| y_{n-2}-y_{n-1}\|^2 
				+ 2 \sum_{n=n_0+1}^{+\infty} \lambda_n\beta_n \left[\sup_{u\in C}\varphi_B \left(u,\frac{p}{\beta_n}\right)-\sigma_C\left(\frac{p}{\beta_n}\right) \right] \\
				&\quad + 2 \left\langle \sum_{n=0}^{K} \lambda_n u - \sum_{n=0}^{K} \lambda_n y_n - \sum_{n=0}^{n_0}\lambda_n u + \sum_{n=0}^{n_0} \lambda_n y_n, w \right\rangle \\
			&= R + 2 \left\langle\sum_{n=0}^{K} \lambda_n u -\sum_{n=0}^{K} \lambda_n y_n -\sum_{n=0}^{n_0} \lambda_n u + \sum_{n=0}^{n_0} \lambda_n y_n, w\right\rangle,
		\end{align*}
	\end{small}
	\begin{fleqn}[\parindent]
		\noindent \begin{align*}
			 \mbox{where}\;  R = \sum_{n=n_0+1}^{+\infty}M_{n-1}^2 \| y_{n-2}-y_{n-1}\|^2   + 2\sum_{n\geq n_0+1}^{+\infty} \lambda_n\beta_n \left[\displaystyle \sup_{u\in C}\varphi \left(u,\frac{p}{\beta_n}\right)-\sigma_C\left(\frac{p}{\beta_n}\right) \right] \in \Rbb.
		\end{align*}
	\end{fleqn}
	Discarding the nonnegative term $\|x_{K+1}-u\|^2$ and dividing $2\tau_{K} = 2\sum_{n=1}^{K}\lambda_n$ we obtain
	\begin{align*}
		-\frac{\|x_{n_0+1}-u\|^2}{2\tau_{K}} 
		&\leq \frac{R+2\inner{-\sum_{n=1}^{n_0} \lambda_n u + \sum_{n=1}^{n_0} \lambda_n y_n}{w}}{2\tau_{K}}	+\frac{2 \inner{\sum_{n=1}^{K} \lambda_n u - \sum_{n=1}^{K} \lambda_n y_n}{w} }{2 \sum_{n=1}^{K}\lambda_n} \\
		&=\frac{\tilde{R}}{2\tau_K} + \inner{u-z'_K}{w},
	\end{align*}
	where $\tilde{R} := R + 2 \inner{-\sum_{n=1}^{n_0} \lambda_n u +\sum_{n=1}^{n_0} \lambda_n x}{w}\in \Rbb$. By passing to the limit as $K\rightarrow +\infty$ and using that $\lim_{K\rightarrow +\infty} \tau_{K} = +\infty$, we get
	\begin{align*}
		\liminf_{K \rightarrow +\infty} \inner{u-z'_K}{w} \geq 0.
	\end{align*}
	Since $z$ is a weak cluster point of $(z'_n)_{n\in\Nbb\cup\{0\}}$, we obtain that $\inner{u-z}{w}\geq 0$. Finally as this inequality holds for arbitrary $(u,w)\in Gr(A+D+N_C)$, $z$ lies in $zer(A+D+N_C)$ (since $A+D+N_C$ is maximally monotone, see preliminaries).
	\\
	
	\noindent (c) Let $n\in \Nbb\cup\{0\}$. It is enough to prove that $\lim_{n\rightarrow +\infty} \|z_n-z'_n\|=0$ and the statement of the theorem will be consequence of Lemma \ref{Lemma1}. Taking $u\in zer(A+D+N_C)$ and $w=0=v+p+D(u)$ where $v\in A(u)$ and $p\in N_C(u)$, it follows from the proof of (a) and the proof of Lemma \ref{Lemma decreasing ineqality} (see (\ref{diff_x_n+1_&_y_n})) that $\sum_{n\in\mathbb{N}\cup\{0\}}  \|y_{n-1}-y_n\|^2 < +\infty$ and $ \| x_{n+1} - y_n \| \leq \left(\frac{\lambda_n\beta_n}{\mu}+\frac{\lambda_n}{\eta} \right) \| y_{n-1}-y_n\| $, respectively. For any $n\in\Nbb\cup\{0\}$ it holds
	\begin{align*}
	\|z_n-z'_n\|^2 
	&= \frac{1}{\tau^2_n} \left\|\sum_{k=0}^{n} \lambda_k (x_k-y_k) \right\|^2\\
	&\leq \frac{1}{\tau_n^2} \left( \sum_{k=0}^{n} \lambda_k \|x_k - y_k\| \right)^2\\
	&\leq \frac{1}{\tau_n^2} \left(\sum_{k=0}^{n} \lambda_k^2 \right) \left( \sum_{k=0}^{n} \|x_k -y_k\|^2\right) \\
	&\leq  \frac{1}{\tau_n^2} \left(\sum_{k=0}^{n} \lambda_k^2 \right) \left( \sum_{k=0}^{n} (2\|x_k -y_{k-1}\|^2 + 2 \|y_{k-1} - y_k\|^2 )\right)\\
	&=   \frac{1}{\tau_n^2} \left(\sum_{k=0}^{n} \lambda_k^2 \right) \left(  2 \sum_{k=0}^{n}\|x_k -y_{k-1}\|^2 + 2 \sum_{k=0}^{n} \|y_{k-1} - y_k\|^2 \right)\\
	&\leq \frac{1}{\tau_n^2} \left(\sum_{k=0}^{n} \lambda_k^2 \right) \left(  2 \sum_{k=0}^{n} \left(\frac{\lambda_{k-1}\beta_{k-1}}{\mu}+\frac{\lambda_{k-1}}{\eta} \right)^2 \| y_{k-2}-y_{k-1}\|^2 + 2 \sum_{k=0}^{n} \|y_{k-1} - y_k\|^2 \right).
	\end{align*}
	Since $\limsup\limits_{n\rightarrow+\infty} \left(\frac{\lambda_n\beta_n}{
		\mu}+\frac{\lambda_n}{
		\mu} \right) < \frac{1}{2}$ and $(\lambda_n)_{n\in \Nbb \cup\{0\}} \in \ell^2 \setminus \ell^1$, taking into consideration that $\tau_n = \sum_{k=0}^{n} \lambda_k \rightarrow + \infty \; (n\rightarrow +\infty) $, obtain $\|z_n-z'_n\| \rightarrow 0 \; (n\rightarrow 0)$. 
\end{proof}


As observed in the forward-backward-forward penalty scheme discussed in [\cite{2014BotCsetnek,BotCsetnek2014}], strong monotonicity of the operator $A$ guarantees the strong convergence of the sequence $(x_n)_{n\in\mathbb{N}\cup\{0\}}$. However, it remains to be seen whether this result extends to the forward-backward-forward with extrapolation from the past in the penalty scheme under consideration. The following result provides the necessary clarification.

\begin{theorem}\label{strongly monotone and strong convergnece}
		Let $(x_n)_{n\in\mathbb{N}\cup \{0\}}$ and $(y_n)_{n\in\mathbb{N} \cup \{0,-1\}}$ be the sequences generated by \textbf{Algorithm 1}. If $(H_{fitz})$ is fulfilled, $\limsup_{n\rightarrow+\infty} \left(\frac{\lambda_n \beta_n}{\mu} + \frac{\lambda_n}{\eta}\right) < \frac{1}{2} $ and the operator $A$ is $\gamma$-strongly monotone with $\gamma > 0$, then $(x_n)_{n\in\mathbb{N}\cup \{0\}}$ converges strongly to the unique element in $zer (A+D+N_{C})$ as $n\rightarrow +\infty$.
	\end{theorem}

	\begin{proof}

 	 Let $u$ be the unique element of $zer (A+D+N_C)$ and $w=0=v+p+D(u)$, where $v\in A(u)$ and $p\in N_C(u)$. We can follow the proof of Lemma \ref{Lemma decreasing ineqality} under $A$ is $\gamma$-strongly monotone with $\gamma > 0$. This means that we obtain $\frac{x_n-y_n}{\lambda_n} - D(y_{n-1}) - \beta_n B(y_{n-1}) \in A(y_n)$ (see (\ref{1_equality of lemma proof})), and one simply show that 
 	 for each $n\in\Nbb$
 	 \begin{align*}
 	 	&2\gamma\lambda_n\|y_n-u\|^2 + \|x_{n+1}-u\|^2-\|x_n-u\|^2 \ + \left( \frac{1}{2} - \left( \frac{\lambda_n\beta_n}{\mu} + \frac{\lambda_n}{\eta} \right)^2  \right) \| y_{n-1}-y_n\|^2 \\
 	 	&\leq \left( \frac{\lambda_{n-1}\beta_{n-1}}{\mu} + \frac{\lambda_{n-1}}{\eta} \right)^2 \| y_{n-2}-y_{n-1}\|^2 +2\lambda_n\beta_n \left[\sup_{u\in C} \varphi_{B} \left(u,\frac{p}{\beta_n}\right)-\sigma_{C} \left(\frac{p}{\beta_n}\right)\right] ,
 	 \end{align*}
 	 equivalently,
 	\begin{align*}
 		&2\gamma\lambda_n\|y_n-u\|^2 + \|x_{n+1}-u\|^2-\|x_n-u\|^2 \ + \left( \frac{1}{2} - 2\left( \frac{\lambda_n\beta_n}{\mu} + \frac{\lambda_n}{\eta} \right)^2  \right) \| y_{n-1}-y_n\|^2 \\
 		&+  \left( \frac{\lambda_n\beta_n}{\mu} + \frac{\lambda_n}{\eta} \right)^2 \|y_{n-1}-y_n\|^2 \\ 
 		&\leq \left( \frac{\lambda_{n-1}\beta_{n-1}}{\mu} + \frac{\lambda_{n-1}}{\eta} \right)^2 \| y_{n-2} -y_{n-1}\|^2 +2\lambda_n\beta_n \left[\sup_{u\in C} \varphi_{B} \left(u,\frac{p}{\beta_n}\right)-\sigma_{C}(\frac{p}{\beta_n})\right]. 
 	\end{align*}
	Let $\bar{E_n} = \|x_{n+1}-u\|^2 + \left( \frac{\lambda_n\beta_n}{\mu} + \frac{\lambda_n}{\eta} \right)^2 \|y_{n-1}-y_n\|^2 < +\infty, \quad \forall n\in\mathbb{N}$. Then:
	
	\begin{align*}
	2\gamma \lambda_n \|y_n-u\|^2 + \bar{E}_n + \left( \frac{1}{2} - 2 \left( \frac{\lambda_n \beta_n}{\mu} + \frac{\lambda_n}{\eta} \right)^2 \right)   \| y_{n-1} - y_n \|^2 \\
	\leq \bar{E}_{n-1} + 2\lambda_n\beta_n \left[\sup_{u\in C} \varphi_{B} \left(u,\frac{p}{\beta_n}\right)-\sigma_{C}\left(\frac{p}{\beta_n}\right)\right].
	\end{align*}
	The hypothesis $\left( \limsup\limits_{n\rightarrow +\infty} \left( \frac{\lambda_n \beta_n}{\mu} + \frac{\lambda_n}{\eta} \right) < \frac{1}{2} \right)$ implies the existence of $n_0 \in \Nbb$ such that for every $n\geq  n_0$
	\begin{align*}
	2 \gamma \lambda_n \|y_n - u\|^2 + \bar{E}_n - \bar{E}_{n-1} \leq  2\lambda_n\beta_n \left[\sup_{u\in C} \varphi_{B}\left(u,\frac{p}{\beta_n}\right)-\sigma_{C}\left(\frac{p}{\beta_n}\right)\right].
	\end{align*}
	Then, the inequality can used for finite summation as: for some $N\geq n_0$ 
	\begin{align*}
		2 \gamma \sum_{n= n_0}^{N} \lambda_n \|y_n - u\|^2 +  \bar{E}_{N}  \leq  \bar{E}_{n_0 -1}   + 2 \sum_{n= n_0}^{N} \lambda_n\beta_n \left[\sup_{u\in C} \varphi_{B}\left(u,\frac{p}{\beta_n}\right)-\sigma_{C}\left(\frac{p}{\beta_n}\right)\right].
	\end{align*}
	Since $\bar{E}_{N} \geq 0$, then we can omit this term on the left-hand side of above inequality. Hence, we derive the following inequality:
	\begin{align}\label{Eq for strong convergence}
	2 \gamma \sum_{n\geq n_0}\lambda_n \|y_n - u\|^2 \leq  \bar{E}_{n_0 -1}   + 2 \sum_{n\geq n_0} \lambda_n\beta_n \left[\sup_{u\in C} \varphi_{B}\left(u,\frac{p}{\beta_n}\right)-\sigma_{C}\left(\frac{p}{\beta_n}\right)\right] < +\infty.
	\end{align}
	It follows that 
	\begin{align}\label{lambda|y-u|^2}
	\sum_{n\in\Nbb} \lambda_n \|y_n - u\|^2 < +\infty.
	\end{align}
	From (\ref{diff_x_n+1_&_y_n}), we can derive
	\begin{align*}
	\|x_n-y_n\|^2 
	&\leq  2 \| x_n -y_{n-1}\|^2 + 2 \|y_n-y_{n-1}\|^2 \\
	&\leq 2  \left( \frac{\lambda_{n-1}\beta_{n-1}}{\mu} + \frac{\lambda_{n-1}}{\eta} \right)^2 \|y_{n-2} - y_{n-1}\|^2 + 2 \|y_n-y_{n-1}\|^2.
	\end{align*}
	Since we have already known from Theorem \ref{Thm_weak (ergodic) convergence} $(a)$ that $\sum_{n\in\mathbb{N}} \| y_{n} - y_{n-1} \|^2 < +\infty$ and $\left( \frac{\lambda_{n-1}\beta_{n-1}}{\mu} + \frac{\lambda_{n-1}}{\eta} \right)^2$ is bounded from above by the hypothesis, then 
 
	\begin{align}\label{|xn-yn|^2}
	\sum_{n\in\mathbb{N}}\|x_n-y_n\|^2 
	&\leq 2 \sum_{n\in\mathbb{N}} \left( \frac{\lambda_{n-1}\beta_{n-1}}{\mu} + \frac{\lambda_{n-1}}{\eta} \right)^2  \|y_{n-2} - y_{n-1}\|^2 + 2 \sum_{n\in\mathbb{N}}\|y_n-y_{n-1}\|^2 < +\infty.
	\end{align}
	It follows from (\ref{lambda|y-u|^2}), (\ref{|xn-yn|^2}), and $\lambda_n$ is bounded from above that
		\begin{align*}
			\sum_{n\in\mathbb{N}} \lambda_n\|x_n-u\|^2  
			\leq 2 \sum_{n\in\mathbb{N}}\lambda_n\|x_n-y_n\|^2 + 2 \sum_{n\in\mathbb{N}}\lambda_n\|y_n - u\|^2 < +\infty.
		\end{align*}	
	Because $\sum_{ { n\in\mathbb{N} } }\lambda_n=+\infty$ and $\lim\limits_{n\rightarrow + \infty} \| x_n - u\| $  exists (proof of Theorem \ref{Thm_weak (ergodic) convergence}  $(a)$), we can conclude that $\lim\limits_{n\rightarrow + \infty} \| x_n - u\| = 0 $.
		

	
	\end{proof}


\section{Minimax problems}\label{Minimax optimization problem}

The minimax optimization problem is a classical problem, where the goal is to find a saddle point for a function in two variables. This problem arises in many applications, including game theory, control theory, and robust optimization. In recent years, there has been a growing interest in the study of minimax problems due to their importance in machine learning and data analysis. In this section, we will consider a class of minimax problems with a convex-concave structure, which can be solved using Tseng's forward-backward-forward method with extrapolation from the past and penalty scheme. We will present the problem formulation, discuss the properties of the objective function, and introduce some numerical methods for solving the problem.

Let $f:\mathcal{H}\times\mathcal{G} \rightarrow \mathbb{R}$ be convex-concave and differentiable, where $f(\cdot, y)$ is convex for all $y\in\mathcal{G}$ and $f(x, \cdot)$  is concave for all $x\in\mathcal{H}$. The \textit{Min-Max problem} (or minimax problem) is the problem in the form:
	\begin{align}\label{Min-Max Prob}
	\min_{x\in\mathcal{H}} \max_{y\in\mathcal{G}} f(x,y).
	\end{align}
	A saddle point of (\ref{Min-Max Prob}) is a vector $(\bar{x}, \bar{y})$ fulfilling:
	\begin{align}\label{saddlepoint def}
	f(\bar{x}, y) \leq f(\bar{x}, \bar{y}) \leq f(x,\bar{y}) \quad \forall x\in\mathcal{H},\; \forall y \in \mathcal{G} .
	\end{align}
	
\noindent Let $F:\mathcal{H}\times \mathcal{G} \rightarrow \mathcal{H}\times \mathcal{G}$ defined by $\begin{pmatrix}  x \\ y  \end{pmatrix} \rightarrow \begin{pmatrix}   \nabla_1 f(x,y) \\ -\nabla_2 f(x,y)\end{pmatrix}$, where $\nabla_1$, $\nabla_2$ are the derivative of $f$ with respect to the first and second component, respectively. Then, writing the optimal conditions for the inequalities in (\ref{saddlepoint def}), we get 
	\begin{align}
	\begin{pmatrix}
	0 \\
	0
	\end{pmatrix}
	=
	\begin{pmatrix}
	\nabla_1 f(\bar{x},\bar{y}) \\ 
	-\nabla_2 f(\bar{x},\bar{y})
	\end{pmatrix},
	\quad \text{i.e.,} \quad	
	(\bar{x}, \bar{y}) \in zer(F).
	\end{align}

	Moreover, $F$ is monotone and also is Lipschitz continuous (if we impose Lipschitz continuous properties on $\nabla_1 f, \nabla_2 f$ ), see [\cite{BC-Book}, Proposition 17.10] and [\cite{rockafellar1997convex}, Theorem 35.1]. In constrained case, the problem in (\ref{Min-Max Prob}) can be written as 
	\begin{align}\label{saddlepoint problem}
	\min_{x\in \mathcal{X}}\max_{y\in \mathcal{Y}} f(x,y),
	\end{align}
	where $\mathcal{X}\subset \mathcal{H}, \mathcal{Y}\subset\mathcal{G}$ are nonempty closed convex sets and a saddle point $(\bar{x}, \bar{y}) \in \mathcal{X} \times \mathcal{Y}$ satisfies:
	\begin{align}
	f(\bar{x}, y) \leq f(\bar{x}, \bar{y}) \leq f(x,\bar{y}) \quad \forall x\in\mathcal{X},\; \forall y \in \mathcal{Y}.
	\end{align}
	The characterization of saddle points of (\ref{saddlepoint problem}) follows from the considerations:
	\begin{align*}
	f(x,\bar{y}) \geq f(\bar{x},\bar{y}) \; \forall x\in\mathcal{X}
 &\Leftrightarrow 0\in\partial (f(\cdot,\bar{y})+ \delta_{\mathcal{X}})(\bar{x}) = \nabla_{1} f(\bar{x},\bar{y})+ N_{\mathcal{X}}(\bar{x}),\\
 	f(\bar{x},y) \leq f(\bar{x},\bar{y}) \; \forall y\in\mathcal{Y}
 &\Leftrightarrow 0\in\partial (-f(\bar{x},\cdot)+ \delta_{\mathcal{Y}})(\bar{y}) = -\nabla_{2} f(\bar{x},\bar{y})+ N_{\mathcal{Y}}(\bar{y}),\\
 	\end{align*}
 	equivalently,
 	\begin{align}
 	\begin{pmatrix}
 	0\\
 	0
 	\end{pmatrix}
 	\in 
 	\begin{pmatrix}
 	N_{\mathcal{X}}(\bar{x}) \\
 	N_{\mathcal{Y}}(\bar{y})
 	\end{pmatrix}
 	+
 	\begin{pmatrix}
	\nabla_{1} f(\bar{x},\bar{y}) \\ 
	-\nabla_{2} f(\bar{x},\bar{y})
	\end{pmatrix}
	=
	\underbrace { N_{\mathcal{X}\times\mathcal{Y}}(\bar{x},\bar{y})}_{\text{maximally monotone}} + \underbrace{F(\bar{x},\bar{y})}_{\text{monotone (Lipschitzian)}}.
 	\end{align}
 	
 	Further we consider more involved minimax problems with constrained sets described by linear equations:
 	\begin{align} \label{constrained minimax problem}
	\min_{ \substack{x\in \mathcal{X} \\ K_1 x = b } } \max_{ \substack{y\in \mathcal{Y} \\ K_2y = b'} } f(x,y) \; \Leftrightarrow \; \min_{ \{ x\in\mathcal{X} : K_1x=b \} } \max_{ \{ y\in\mathcal{Y} : K_2y = b' \} } f(x,y),
 	\end{align}
	where $\mathcal{Q}$, $\mathcal{Q}'$ are other real Hilbert spaces, $K_1: \mathcal{H}\rightarrow\mathcal{Q}, K_2:\mathcal{G}\rightarrow\mathcal{Q}'$ are linear continuous operators with closed ranges such that $b\in\mathcal{Q}, b'\in \mathcal{Q}'$ and $\{x\in\mathcal{X} : K_1x = b\}\neq \emptyset$, $\{y\in\mathcal{Y} : K_2 y = b'\}\neq \emptyset$ are nonempty sets. This kind of the minimax constrained problem was studied with a plenty of applications, for instance in [\cite{Yang2022LiLan}]. In this case we have the following characteristic of saddle points of (\ref{constrained minimax problem})
	\begin{align}
	\begin{pmatrix}
 	0\\
 	0
 	\end{pmatrix}
 	\in 
 	\begin{pmatrix}
 	N_{\mathcal{X}\cap \{ x\in\mathcal{X} : K_1x=b \} }(\bar{x}) \\
 	N_{\mathcal{Y}\cap  \{ y\in\mathcal{Y} : K_2y = b' \} }(\bar{y}) 
 	\end{pmatrix}+ 
	\underbrace{ 	
 	\begin{pmatrix}
	\nabla_{1} f(\bar{x},\bar{y}) \\ 
	-\nabla_{2} f(\bar{x},\bar{y})
	\end{pmatrix}
	}_{:=F(\bar{x},\bar{y})}.
	\end{align}
	Note that $N_{C_1\cap C_2} = \partial (\delta_{C_1} + \delta_{C_2}) = N_{C_1} + N_{C_2}$, if a constraint qualification is satisfied, for example $ \emptyset \neq C_1 \cap \operatorname{int} C_2 $ (see [\cite{BC-Book}, Corollary 16.38]). Then
	\begin{align}\label{split monotone inclusion for convex opt prob}
	\begin{pmatrix}
 	0\\
 	0
 	\end{pmatrix}
 	\in 
 	\underbrace{ 
 	\begin{pmatrix}
 	N_{\mathcal{X} }(\bar{x}) \\
 	N_{\mathcal{Y} }(\bar{y})
 	\end{pmatrix} 
 	}_{\text{maximally monotone}}
 	+
 	\begin{pmatrix}
	N_{ \{ x : K_1x=b \} }(\bar{x}) \\
 	N_{ \{ y : K_2y = b' \} }(\bar{y})
 	\end{pmatrix}
 	 + \underbrace{\begin{pmatrix}
	\nabla_{1} f(\bar{x},\bar{y}) \\ 
	-\nabla_{2} f(\bar{x},\bar{y})
	\end{pmatrix}
	}_{\substack{:= F(\bar{x},\bar{y}), \\ \text{monotone (Lipschitzian)} } }.	 	
	\end{align}
	We find the solution of \begin{align}\label{linear equation for constrant minmax}
	K_1x = b, \quad K_2y = b',
	\end{align} by considering the functions $\psi_1:\mathcal{H} \rightarrow \mathbb{R}, \psi_2 : \mathcal{G}\rightarrow\mathbb{R}$ as 
	\begin{align}
	\psi_1(x) = \frac{1}{2}\|K_1x - b\|^2\quad \mbox{and} \quad  \psi_2(y) = \frac{1}{2}\|K_2 y - b'\|^2. 
	\end{align}
	Namely, the solution of (\ref{linear equation for constrant minmax}) can be obtained by solving the following optimization problem:
	\begin{align*}
	\arg\min_{x\in\mathcal{H}} \frac{1}{2}\|K_1x - b\|^2 = \arg\min_{x\in\mathcal{H}} \psi_1(x) \Leftrightarrow \nabla\psi_{1}(x) = 0,\\
	\arg\min_{y\in\mathcal{G}} \frac{1}{2}\|K_2 y - b'\|^2 = \arg\min_{y\in\mathcal{G}} \psi_2(y) \Leftrightarrow \nabla\psi_{2}(y)=0.
	\end{align*}	
	For every $x\in\mathcal{H}$ and $y\in\mathcal{G}$, let 
	\begin{align} \label{setting operators in minimax}
	\bar{A}(x,y) = 
		\begin{pmatrix}
			N_{\mathcal{X}} (x) \\
			N_{\mathcal{Y}} (y)
		\end{pmatrix},
	\;
	\bar{B}(x,y)= 
		\begin{pmatrix}
			\nabla \psi_1 (x) \\
			-\nabla \psi_2 (y)
		\end{pmatrix},
	\;
	\bar{C}= zer \bar{B},
	\;
	\bar{D} (x,y) = 
		\begin{pmatrix}
		 \nabla_{1} f(x,y) \\
		 -\nabla_{2} f(x,y) \\
		\end{pmatrix}	 
	\end{align}

	Therefore, (\ref{split monotone inclusion for convex opt prob}) can be written as  
	\begin{align} \label{monotone inclusion of convex optimization}
	\begin{pmatrix}
 	0\\
 	0
 	\end{pmatrix}
 	\in 
 	\begin{pmatrix}
 	N_{\mathcal{X} }(\bar{x}) \\
 	N_{\mathcal{Y} }(\bar{y})
 	\end{pmatrix} 
 	+ 	
 	\begin{pmatrix}
	N_{zer(\nabla \psi_1)} (\bar{x}) \\
 	N_{zer(\nabla \psi_2)}(\bar{y})
 	\end{pmatrix} 	
 	+ \begin{pmatrix}
	\nabla_{1} f(\bar{x},\bar{y}) \\ 
	-\nabla_{2} f(\bar{x},\bar{y})
	\end{pmatrix} 
	\Leftrightarrow 
	\begin{pmatrix}
 	0\\
 	0
 	\end{pmatrix} \in \bar{A}(\bar{x},\bar{y}) + N_{\bar{C}}(\bar{x},\bar{y}) + \bar{D}(\bar{x},\bar{y}).
	\end{align}
	Due to the above setting, the problem (\ref{split monotone inclusion for convex opt prob}) turns into the monotone inclusion in (\ref{monotone inclusion problem_1}). The iterative pattern in \textbf{Algorithm 1} can be transformed into for each $n\geq 0$ as
	\begin{align*}
	\begin{pmatrix}
	r_n \\
	s_n
	\end{pmatrix}
	&= J_{\lambda_n {\tiny \begin{pmatrix} N_{\mathcal{X}} \\ N_{\mathcal{Y}}  \end{pmatrix} } }
	\left[
		\begin{pmatrix}
		x_n \\ y_n
		\end{pmatrix}
		-
		\lambda_n
		\begin{pmatrix}
		\nabla_{1} f(r_{n-1}, s_{n-1})\\
		-\nabla_{2} f(s_{n-1}, s_{n-1})
		\end{pmatrix}
		-
		\lambda_n \beta_n 
		\begin{pmatrix}
		\nabla \psi_1(r_{n-1}) \\
		\nabla \psi_2(s_{n-1})
		\end{pmatrix}
	\right];\\
	\begin{pmatrix}
	x_{n+1} \\
	y_{n+1}
	\end{pmatrix}
	&=
	\lambda_{n} \beta_{n} 
	\left[
		\begin{pmatrix}
			\nabla\psi_1(r_{n-1}) \\
			\nabla\psi_2(s_{n-1})
		\end{pmatrix}
		-
		\begin{pmatrix}
			\nabla\psi_1(r_n)\\
			\nabla\psi_2(s_n)
		\end{pmatrix}
	\right]
	+
	\lambda_n
	\left[
		\begin{pmatrix}
			\nabla_{1} f(r_{n-1},s_{n-1}) \\
			-\nabla_{2} f({r_{n-1}},s_{n-1} )
		\end{pmatrix}
		-
		\begin{pmatrix}
			\nabla_{1} f(r_{n},s_{n}) \\
			-\nabla_{2} f({r_{n}},s_{n} )
		\end{pmatrix}
	\right]
	+
	\begin{pmatrix}
	r_n \\
	s_n
	\end{pmatrix},
	\end{align*}
 	where $\nabla \psi_1 (x) = K_1^{*} (K_1x-b)$ and $\nabla \psi_2 (y) = K_2^{*} (K_2 y-b')$. 
	Notice that if $G$ is closed convex set subset of $\mathcal{H}$, then $J_{N_{G}} = (Id+N_{G})^{-1} = prox_{\delta_{G}} = proj_G$ (see [\cite{BC-Book}, Example 23.4]). Therefore we are able to construct a relevant algorithm for solving the monotone inclusion scheme in (\ref{monotone inclusion of convex optimization}) (i.e., for solving problem (\ref{constrained minimax problem})) demonstrated as:
	
	\begin{tcolorbox}[colback=white]
		\textbf{Algorithm 2}:
		\begin{equation}
		\begin{aligned} 
			\mbox{Initialization}&:\; \mbox{Choose}\; r_{-1}, x_0 \in\mathbb{R}^{m_1}, s_{-1},y_{0} \in \mathbb{R}^{m_2} \;(\mbox{with}\; r_{-1}= x_0, s_{-1}=y_{0}).   \\
			\mbox{For}\;  n\in\mathbb{N}\cup\{0\}\; \mbox{set}
			&:\; r_n = proj_{\mathcal{X}} [ x_n-\lambda_n \nabla_{1} f(r_{n-1},s_{n-1}) -  \lambda_n\beta_n  K_1^{*} (K_1r_{n-1}-b) ];\\
			&\quad s_n = proj_{\mathcal{Y}} [ y_n + \lambda_n  \nabla_{2} f( r_{n-1}, s_{n-1}) - \lambda_n\beta_n   K_1^{*} (K_2s_{n-1}-b') ];\\
			&\quad x_{n+1}= \lambda_n\beta_n [ K_1^{*} (K_1r_{n-1}-b) - K_1^{*} (K_1r_{n}-b)  ] \!+\! \lambda_n  \nabla_{1} f(r_{n-1},s_{n-1}) \\ 			&\quad\quad\quad\quad -\lambda_n  \nabla_{1} f(r_n,s_n)  + r_n;\\
			&\quad y_{n+1}= \lambda_n\beta_n [  K_1^{*} (K_2s_{n-1}-b') - K_1^{*} (K_2s_{n}-b') ] \!+\! \lambda_n \nabla_{2} f(r_n,s_{n});\\
			&\quad\quad\quad\quad - \lambda_n  \nabla_{2} f(r_{n-1},s_{n-1})+s_n,
		\end{aligned}
		\end{equation}
	\end{tcolorbox}
	\noindent where  $(\lambda_n)_{n\in\mathbb{N}\cup\{0\}}$ and  $(\beta_n)_{ { n\in\mathbb{N}\cup\{0\} } }$  are positive real numbers.
	
	We show that we can use condition $(ii)$ in $(H)$ on page 1 in order to get the convergent results. We define $\Psi : \mathcal{H} \times \mathcal{G} \rightarrow \mathbb{R}$,
		\begin{align}
			\Psi(x,y) := \psi_1(x) + \psi_2 (y) = \frac{1}{2}\| K_1 x- b\|^2 +  \frac{1}{2}\| K_2 - b' \|^2.
		\end{align}
	 Then, we have (see [\cite{BC-Book}, Proposition 13.27])
	 	\begin{align}
	 		\Psi^* (x, y) = \psi_1^*(x) + \psi_2^*(y).
 		\end{align}
 	From the considerations above we have:
		\begin{align*}
			\bar{C} = \left\{ (x,y) : 
					\begin{pmatrix}
					\nabla \psi_1(x) \\
					-\nabla \psi_2(y)
					\end{pmatrix}
				=
					\begin{pmatrix}
					0\\
					0
					\end{pmatrix}
			\right\}
			= \left\{ (x,y) : 
			\begin{matrix}
				K_1 x = b \\
				K_2 y = b'
			\end{matrix} \right\}
		\;
		\mbox{and}
		\;
		 \bar{C} = zer \psi_1 \times zer \psi_2 = zer \Psi. 
		\end{align*}
 	Furthermore, we also have
 		\begin{align*}
 			\sigma_{\bar{C}} = \sigma_{zer \psi_1 \times zer \psi_2} = \sigma_{zer \psi_1} + \sigma_{zer \psi_2} \quad \mbox{and} \quad N_{\bar{C}} = N_{zer{\psi_1}} \times N_{zer\psi_2}.
 		\end{align*}
 	Thus, if we impose for $i=1,2$
 		\begin{align}\label{bounded for i in minimax problem}
 			\forall p_i \in \operatorname{ran} (N_{zer \psi_i}), \quad \sum_{n\in\mathbb{N}\cup\{0\}} \lambda_{n}\beta_{n} \left[\psi_i^* (\frac{p_i}{\beta_n}) - \sigma_{zer \psi_i} (\frac{p_i}{\beta_{n}}) \right] < + \infty,
 		\end{align}
 	then, it follows that
 		\begin{align*}
 			\forall (p_1,p_2)  \in \operatorname{ran} (N_{zer \bar{C}}), \quad \sum_{n\in\mathbb{N}\cup\{0\}} \lambda_{n}\beta_{n} \left[ \Psi^* \left(\frac{(p_1,p_2)}{\beta_n}\right) - \sigma_{\bar{C}} \left(\frac{(p_1,p_2)}{\beta_{n}}\right)  \right]< + \infty,
 		\end{align*}
	exactly the condition (ii) on page 1 that we need for convergence in $(H)$. 
	Hence, we can propose the ergodic convergence results which follow from Thorem \ref{Thm_weak (ergodic) convergence} as below.
	\begin{corollary}
		In the framework of (\ref{constrained minimax problem}), let $\mathcal{X}\subset \mathcal{H}, \mathcal{Y}\subset\mathcal{G}$ be nonempty closed convex sets and the operators $\bar{A}, \bar{B}, \bar{C}$ and $\bar{D}$  given as in (\ref{setting operators in minimax}) where $\bar{A}$ is maximally monotone, $\bar{D}$ is monotone and $\eta^{-1}$-Lipschitz continuous with $\eta > 0$, $\bar{B}$ is monotone and $\mu^{-1}$-Lipschitz continuous operator with $\mu > 0$ (i.e., $\mu^{-1} :=  \bar{K}$, where $\bar{K} = \max\{\| K_1\|^2, \| K_2\|^2\}, K_1\neq 0, K_2 \neq 0$), and $\bar{C} :=zer \bar{B} \neq \emptyset$. Consider the sequences generated by \textbf{Algorithm 2}, and the sequence $(z_n)_{n\in\mathbb{N}\cup\{0\}}$ defined as 
		\begin{align*}
				z_n = \frac{1}{\tau_n} \sum_{k=0}^{n} \lambda_k \begin{pmatrix}
					x_k \\
					y_k
				\end{pmatrix},\quad \mbox{where}\; \tau_n = \sum_{k=0}^{n} \lambda_k \; \forall n \in \mathbb{N},
		\end{align*} 
		where $(\lambda_{k})_{k\in \mathbb{N}\cup\{0\}}$ is a positive sequence such that $\sum_{k\in\mathbb{N}\cup\{0\}} \lambda_{k} = +\infty$. If we suppose the following hypotheses:
		\begin{align*}
			(H^{min-max})\begin{cases}
					(i) &\bar{A}+ N_{\bar{C}}  \;\mbox{is maximally monotone and} \; (\ref{monotone inclusion of convex optimization}) \; \mbox{has a solution};\\
					(ii) &\mbox{For every}\;  p_i \in \operatorname{ran} (N_{zer \psi_i}), \sum_{n\in\mathbb{N}\cup\{0\}} \lambda_{n}\beta_{n} \left[\psi_i^* (\frac{p_i}{\beta_n}) - \sigma_{zer \psi_i} (\frac{p_i}{\beta_{n}}) \right] < + \infty, i= 1,2;\\
					(iii) &(\lambda_n)_{ n\in\mathbb{N}\cup\{0\} }  \in \ell^2 \setminus \ell^1,
				\end{cases}
		\end{align*} and $\limsup_{n\rightarrow +\infty} \left( \frac{\lambda_n \beta_n}{\mu} + \frac{\lambda_n}{\eta} \right)< \frac{1}{2}$, then $(z_n)_{n\in\mathbb{N}\cup\{0\}}$ converges weakly to a saddle point of (\ref{constrained minimax problem}).
	\end{corollary}

	\begin{remark}\
		\begin{enumerate}
			\item Considering the condition $(i)$ in $(H^{min-max})$, notice that $\bar{A}+ N_{\bar{C}}$ is maximally monotone if a regularity condition is fulfilled, for example (see [\cite{BC-Book}, Corollary 24.4])
			\begin{equation}\label{Condition A}\leqnomode
			(\operatorname{int}\mathcal{X}\times \operatorname{int} \mathcal{Y}) \cap \left\{ (x,y) : 
			\begin{matrix}
				K_1 x = b \\
				K_2 y = b'
			\end{matrix} \right\}
			\neq \emptyset.  \tag{Condition A}
			\end{equation}
			In addition, $(\ref{monotone inclusion of convex optimization})$ has a solution if the (\ref{Condition A}) holds and (\ref{constrained minimax problem}) has saddle points.
			
			\item Notice that $(ii)$ in $(H^{min-max})$ is fulfilled when $\sum_{n\in\mathbb{N}\cup\{0\}} \frac{\lambda_n}{\beta_{n}} < + \infty$ and $b=b'=0$, see Section \ref{Motivation}.
			
			\item 	We can show the Lipschitz property of $\bar{B}$ as follows:
			\begin{align*}
				\| \bar{B}(x,y) - \bar{B}(u,v) \| 
				&= \left\|  \begin{pmatrix}  \nabla \psi_1(x) \\ -\nabla \psi_{2}(y)\end{pmatrix} - \begin{pmatrix} \nabla \psi_1(u) \\ - \nabla\psi_{2} (v) \end{pmatrix} \right\| \\
				&= \left\| \begin{pmatrix}  \nabla \psi_1(x) - \nabla \psi_1(u) \\ \nabla \psi_{2}(v)-\nabla \psi_{2}(y)\end{pmatrix} \right\| \\
				&= \left\| \begin{pmatrix}  K_1^* (K_1x-b) - K_1^*(K_1u-b)  \\ K_2^* (K_2v-b') - K_2^*(K_2y-b') \end{pmatrix} \right\| \\
					&= \left\| \begin{pmatrix}   K_1^*K_1\left(x - u \right) \\ - K_2^*K_2\left(y -  v \right)  \end {pmatrix} \right\|\\
						&\leq  \|\bar{K}\|  \|(x,y) - (u,v) \|,
					\end{align*}
					where $\bar{K} = \max\{\| K_1\|^2, \| K_2\|^2\}$.
		\end{enumerate}
	
	\end{remark}

\section{The forward-backward-forward algorithm with extrapolation from the past and penalty scheme for the problem involving composion of linear continuous operators}\label{The FBF algorithm with extrapolation penalty scheme for the problem involving linearly composed and parallel-sum type monotone operators}

In this section, we propose forward-backward-forward algorithm with extrapolation from the past and penalty scheme utilized to address the inclusion problem involving the finite sum of the compositions of monotone operators with linear continuos operators. 
To this end, we begin with the following problem:\\

\noindent \textbf{Problem 3:} Let $\mathcal{H}$ be a real Hilbert space, $A : \mathcal{H}\rightrightarrows \mathcal{H}$ a maximally monotone operator, $D : \mathcal{H}\rightarrow \mathcal{H}$ a monotone and $\nu$-Lipschitz continuous operator with $\nu >0$. Let $m$ be strictly positive integer and for any $i\in\{ 1,\dots,m \}$, let $\mathcal{G}_i$ be a real Hilbert space, $B_i : \mathcal{G}_i \rightrightarrows \mathcal{G}_i$ a maximally monotone operator, 
 and $L_i : \mathcal{H}\rightarrow \mathcal{G}_i$ a nonzero linear continuous operator. Assume that $B:\mathcal{H}\rightarrow \mathcal{H}$ is a monotone and $\mu^{-1}$-Lipschitz continuous operator with $\mu >0$ and suppose $C = zerB\neq\emptyset$. The monotone inculsion problem to solve is
	\begin{align} \label{monotone inclusion involing composition of linear continuous operators}
	0 \in Ax + \sum_{i=1}^{m} L_i^*B_i 
	 L_i x + Dx + N_C(x).
	\end{align}
	
	The monotone inclusion problem, \textbf{Problem 3}, can be reformulated into the same form as \textbf{Problem 1} using the product space approach with a pertinent setting. To address this problem, we work with the product space $\mathcal{H} \times \mathcal{G}_1 \times \dots \times \mathcal{G}_m$, where the inner product and associated norm are defined for every element $(x,v_1,\dots,v_m)$ by $\inner{(x,v_1,\dots,v_m)}{(y,w_1,\dots,w_m)} = \inner{x}{y} + \sum_{i=1}^{m} \inner{v_i}{w_i}$ and $\|(x,v_1,\dots,v_m)\| = \sqrt{\|x \|^2 + \sum_{i=1}^{m}\|v_i \|^2}$. We proceed to define the operators on the product space $\mathcal{H} \times \mathcal{G}_1 \times \dots \times \mathcal{G}_m$ as follows: for every $(x,v_1,\dots,v_m), (y,w_1,\dots,w_m) \in \mathcal{H} \times \mathcal{G}_1 \times \dots \times \mathcal{G}_m$, $\tilde{A}: \mathcal{H} \times \mathcal{G}_1 \times \dots \times \mathcal{G}_m \rightrightarrows \mathcal{H} \times \mathcal{G}_1 \times \dots \times \mathcal{G}_m$ 
	
	\begin{align*}
	\tilde{A}(x,v_1,\dots,v_m) = Ax \times B_{1}^{-1} v_1 \times \dots \times B_{m}^{-1} v_m;
	\end{align*}	 
	$\tilde{D} : \mathcal{H} \times \mathcal{G}_1 \times \dots \times \mathcal{G}_m \rightarrow \mathcal{H} \times \mathcal{G}_1 \times \dots \times \mathcal{G}_m$,
	\begin{align*}
	\tilde{D}(x,v_1,\dots,v_m) = \left( \sum_{i=1}^{m} L_i^* v_i + Dx, 
	- L_1 x, \dots, 
	- L_m x \right);
	\end{align*}
	and $\tilde{B} : \mathcal{H} \times \mathcal{G}_1 \times \dots \times \mathcal{G}_m \rightarrow \mathcal{H} \times \mathcal{G}_1 \times \dots \times \mathcal{G}_m$,
	\begin{align*}
	\tilde{B}(x,v_1,\dots,v_m) = (Bx,\underbrace{0,\dots,0}_{m-\text{elements}}).
	\end{align*}
	
	Note that the maximal monotonicity of $A$ and $B_i$, $i=1,\dots,m$, implies that the operator $\tilde{A}$ is also maximally monotone, as stated in [\cite{BC-Book}, Proposition 20.22 and 20.23]. Moreover, according to  [\cite{CombettesPesquet2012}, Theorem 3.1], it is possible to demonstrate that the operator $\tilde{D}$ is monotone and $\beta$-Lipschitz continuous 
	where $\beta = \nu + \sqrt{\sum_{i=1}^{m} \|L_i\|^2 }$ (in this context, the monotonicity and Lipschitzian of $D$ is a special case of [\cite{2014BotCsetnek,CombettesPesquet2012}]). Moreover, we also have that  $\tilde{B}$ is monotone and $\mu^{-1}$-Lipschitz continuous with
	\begin{align*}
	\tilde{C} := zer \tilde{B} := zer {B} \times \mathcal{G}_1 \times \dots \times \mathcal{G}_m = C \times \mathcal{G}_1 \times \dots \times \mathcal{G}_m ,
	\end{align*}
	and
	\begin{align*}
	N_{\tilde{C}} (x,v_1,\dots,v_m) = N_{C}(x)\times \underbrace{\{0\} \times \dots\times \{0\} }_{m\text{-elements}}.
	\end{align*}

	One can show (see [\cite{2014BotCsetnek}]) that $x$ is a solution to \textbf{Problem 3} if and only if there exists $v_1\in\mathcal{G}_1,\dots,v_m\in\mathcal{G}_m$ such taht $(x,v_1,\dots,v_m) \in zer(\tilde{A} + \tilde{D} + N_{\tilde{C}})$. On the other hand, when $(x,v_1,\dots,v_m) \in zer(\tilde{A} + \tilde{D} + N_{\tilde{C}})$, then $x\in zer(A + \sum_{i=1}^{m} L_i^{*} B_i  L_i + D + N_C)$. This means that determining the zeros of $\tilde{A} + \tilde{D} + N_{\tilde{C}}$ will automatically provide a solution to \textbf{Problem 3}.

	By using the identity given in [\cite{BC-Book}, Proposition 23.16], which provides that for any $(x,v_1,\dots,v_m)\in \mathcal{H} \times \mathcal{G}_1, \times \dots \times \mathcal{G}_m$ and $\lambda > 0$, $J_{\lambda\tilde{A}} (x,v_1,\dots,v_m) = (J_{\lambda_n A}(x), J_{\lambda B_{1}^{-1}} (v_1),\dots,J_{\lambda B_{m}^{-1}} (v_m))$, it can be observed that the iterations of \textbf{Algorithm 1} are given for any $n\in\mathbb{N}\cup\{0\}$ as follows:
	
		\begin{align*}
	(y_n, q_{1,n}, \dots, q_{m,n}) 
	&= J_{\lambda_n \tilde{A}} [ (x_n, v_{1,n},\dots,v_{m,n}) - \lambda_n \tilde{D}(y_{n-1},q_{1,n-1},\dots,q_{m,n-1})\\ 
	&\quad - \lambda_{n}\beta_{n}\tilde{B} (y_{n-1}, q_{1,n-1},\dots,q_{m,n-1})],\\
	(x_{n+1},v_{1,n+1},\dots,v_{m,n+1})&= \lambda_n \beta_n [\tilde{B}(y_{n-1}, q_{1,n-1}, \dots, q_{m,n-1}) - \tilde{B}(y_{n}, q_{1,n}, \dots, q_{m,n})] \\
	&\quad + \lambda_n [\tilde{D}(y_{n-1}, q_{1,n-1}, \dots, q_{m,n-1}) - \tilde{D}(y_{n}, q_{1,n}, \dots, q_{m,n})] \\
	&\quad +  (y_{n}, q_{1,n}, \dots, q_{m,n}). 
	\end{align*}
	
		This leads to the algorithm below:
	
	\begin{tcolorbox}[colback=white]
		\textbf{Algorithm 3}:
		\begin{equation}
			\begin{aligned} 
				\mbox{Initialization}&:\; \mbox{Choose}\;  \;(x_0,v_{1,0},\dots,v_{m,0}), (y_{-1},q_{1,-1},\dots,q_{m,-1}) \in \mathcal{H} \times \mathcal{G}_1,\! \times \! \dots \! \times \! \mathcal{G}_m \\ 
				&\quad \; \mbox{with} \;  y_{-1}=x_{0},q_{i,-1}=v_{i,0} \;\forall i = 1,\dots,m. \\
				\mbox{For}\; n\in\mathbb{N}\cup\{ 0 \}&, \\ 
					\mbox{set} \; &:\; y_n = J_{\lambda_n A}[ x_n-\lambda_n (D(y_{n-1})+ \sum_{i=1}^{m} L^* (q_{i,n-1}) )-\lambda_n\beta_n B(y_{n-1})];\\
					&\quad  q_{i,n} = J_{\lambda_n B_i^{-1}} [v_{i,n} + \lambda_n L_i (y_{n-1})], \; i=1,\dots,m,  ; \\ 
					&\quad x_{n+1}\!=\!\lambda_n\beta_n[B(y_{n-1})\!-\!B(y_n)] \!+\! \lambda_n[D(y_{n-1})\!-\!D(y_n)]\!\\
					&\quad\quad\quad\quad+\! \lambda_n \sum_{i=1}^{m} L^*_i(q_{i,n-1}-q_{i,n})\!+\!y_n;\\
					&\quad v_{i,n+1} \!=\! \lambda_n L_i(y_n-y_{n-1}) 
					\!+\! q_{i,n}, \; i=1,\dots,m,  
			\end{aligned}
		\end{equation}
	\end{tcolorbox}
	\noindent where  $(\lambda_n)_{n\in\mathbb{N}\cup\{0\}}$ and  $(\beta_n)_{ { n\in\mathbb{N}\cup\{0\} } }$  are positive real numbers. To consider the convergence of this iterative scheme, we need the following additionally hypotheses, similar to the hypotheses in [\cite{2014BotCsetnek}, Section 2]:
	
	$(H_{fitz}^{sum})\begin{cases}
		(i) &A +N_C  \;\mbox{is maximally monotone and} \; zer(A+\sum_{i=1}^{m} L^*B_i
		L+D+N_C)\neq\emptyset;\\
		(ii) &\mbox{For every } p\in \operatorname{ran} (N_C),  \sum_{n\in\mathbb{N}\cup\{0\}}  \lambda_n\beta_n [\sup\limits_{u\in C} \varphi_{B} (u,\frac{p}{\beta_n})-\sigma_C(\frac{p}{\beta_n})] < +\infty;\\
		(iii) &(\lambda_n)_{ n\in\mathbb{N}\cup\{0\} }  \in \ell^2 \setminus \ell^1.
	\end{cases}$
	
	The Fitzpatrick function of $\tilde{B}$ (i.e., $\varphi_{\tilde{B}}$ ) and the support function of $\tilde{C}$ (i.e., $\sigma_{\tilde{C}}$) can be computed in the same way as demonstrated in [\cite{2014BotCsetnek}, Section 2] for arbitrary elements \emergencystretch 3em $(x_1,v_{1},\dots,v_{m}), (x_1',v_{1}',\dots,v_{m}')\in \mathcal{H} \times \mathcal{G}_1, \times \dots \times \mathcal{G}_m$. Moreover, the satisfaction of condition $(i)$ in $(H_{fitz}^{sum})$ guarantees that $\tilde{A} + N_{\tilde{C}}$ is maximally monotone and $zer(\tilde{A}+\tilde{D}+N_{\tilde{C}})\neq\emptyset$. As a result, we can apply Theorem \ref{Thm_weak (ergodic) convergence} and \ref{strongly monotone and strong convergnece} to the problem of finding the zeros of $\tilde{A}+\tilde{D}+N_{\tilde{C}}$. Thus, we can establish the convergence results for the scheme of the sum of monotone operators and linear continuous operators, which are given below.

	\begin{theorem}
	Let $(x_n)_{n\in\mathbb{N}\cup\{0\} }$ and $(z_n)_{n\in\mathbb{N}\cup\{0\} }$ be sequences generated by \textbf{Algorithm 3}. Assume $(H_{fitz}^{sum})$ is fulfilled and $\limsup_{n\rightarrow +\infty} (\frac{\lambda_n \beta_n}{\mu} + \lambda_n \beta) < \frac{1}{2}$, where
	$$\beta = \nu
	+ \sqrt{\sum_{i=1}^{m}  \|L_i\|^2}.$$
	Then $(z_n)_{n\in\mathbb{N}\cup\{0\} }$ converges weakly to an element in $zer(A + \sum_{i=1}^{m} L_i^{*} B_i
	 L_i + D + N_C)$ as $n\rightarrow\infty$. If, additionally, $A$ and $B^{-1}_i, i=1,\dots,m$ are strongly monotone, then $(x_n)_{n\in\mathbb{N}\cup\{0\} }$ converges strongly to the unique element in $zer(A + \sum_{i=1}^{m} L_i^{*} B_i
	 L_i + D + N_C)$  as $n\rightarrow\infty$.
	\end{theorem}	
	
	\begin{remark}
		In case $m=1$, our considered problem will become a solver of the monotone inclusion problems involing compositions with linear continous operators suggested in [\cite{BotCsetnek2014}, Section 3.2] corresponding to the context of the iterative method based on Tseng's forward-backward-forward algorithm with extrapolation from the past and   penalty scheme.
	\end{remark}


\section{Convex minimization problem}\label{Convex minimization problem}

In this section, we apply the obtained results by using the forward-backward-forward algorithm with extrapolation from the past and penalty scheme for monotone inclusion problems to the minimization of a convex function with a complex formulation subject to the set of minimizers of another convex and differentiable function with a Lipschitz continuous gradient. We consider the convex minimization problem proposed in [\cite{2014BotCsetnek}], given by:\\

\noindent \textbf{Problem 4:} Let $\mathcal{H}$ be a real Hilbert space, $f\in\Gamma(\mathcal{H})$ and $h:\mathcal{H} \rightarrow \mathbb{R}$ a convex and differentiable function with $\nu$-Lipschitz continous gradient for $\nu>0$. Let $m$ be a strictly positive integer and for any $i = 1,\dots,m$, let $\mathcal{G}_i$ be a real Hilbert space, $g_i\in \Gamma(\mathcal{G}_i)$, 
 and $L_i:\mathcal{H}\rightarrow \mathcal{G}_i$ a nonzero linear continuous operator. Furthermore, let $\Psi: \mathcal{H}\rightarrow \mathbb{R}$ be convex and differentiable with a $\mu^{-1}$-Lipschitz continuous gradient, fulfilling $\min \Psi = 0$. The convex minimization problem under investigation is

	\begin{align}\label{Convex-Opt Promblem}
	\inf_{x\in \arg\min \Psi} \bigg\{ f(x) + \sum_{i=1}^{m} g_i
	 (L_i x) + h(x) \bigg\}.
	\end{align}	 
	
	In order to solve this problem in our context, we can follow the \textbf{Problem 3} in section \ref{The FBF algorithm with extrapolation penalty scheme for the problem involving linearly composed and parallel-sum type monotone operators}. We khow that any element belonging to $zer(\partial f + \sum_{i=1}^{m} L^*_i \partial g_i 
	 L_i + \nabla h + N_C)$ is an optimal solution for (\ref{Convex-Opt Promblem}) if we substitute all of operators with
	\begin{align}\label{setting operators for convex optimization}
		A =\partial f, \; B = \nabla \Psi,\; C = \arg\min \Psi = zer B,\; D = \nabla h, \;\mbox{and}\; B_i = \partial g_i, 
		\; i = 1,\dots,m,
	\end{align}
	where $B$ is a monotone and $\mu^{-1}$-Lipschitz continuous operator (see [\cite{BC-Book}, Proposition 17.10, Theorem 18.15]). However, the converse is true only if a suitable qualification condition is satisfied. Necessary conditions, namely that (\ref{Convex-Opt Promblem}) leads to \textbf{Problem 3} with (\ref{monotone inclusion involing composition of linear continuous operators}), have been previously given in [\cite{2014BotCsetnek}] 
	(for example, see [\cite{CombettesPesquet2012}, Proposition 4.3, Remark 4.4])
	\begin{align}\label{sqri-qualification condition}
		(0,\dots,0) \in \operatorname{sqri} \bigg( \prod^{m}_{i=1} \operatorname{dom} g_i - \{ (L_1x,\dots,L_mx) : x\in \operatorname{dom}\; f \cap C \} \bigg).
	\end{align}	

	Note that the \textit{strong quasi-relative interior} for a convex set $S$ in a real Hilbert space $\mathcal{H}$ is denoted by $\operatorname{sqri} S$ and is defined as the set of points $x\in S$ such that the union of all positive scalar multiples of $(S-x)$ is a closed linear subspace of $\mathcal{H}$. Notice that $\operatorname{sqri} S$ always contains the \textit{interior} of $S$, denoted by $\operatorname{int} S$, but this inclusion can be strict. When $\mathcal{H}$ is finite-dimensional, $\operatorname{sqri} S$ coincides with the \textit{relative interior} of $S$, denoted by $\operatorname{ri} S$, which is the interior of $S$ with respect to its affine hull. Moreover, the condition in (\ref{sqri-qualification condition}) is fulfilled if
	\begin{itemize}
	\item $\operatorname{dom} \; g_i = \mathcal{G}_i, i = 1,\dots,m$ or
	\item $\mathcal{H}$ and $\mathcal{G}_i$ are finite dimensional and there exists $x\in \operatorname{ri} \operatorname{dom}f \cap  \operatorname{ri}  C$ such that $L_i x \in \operatorname{ri} \operatorname{dom} \; g_i, i=1,\dots,m$ (see [\cite{CombettesPesquet2012}, Proposition 4.3]).
	\end{itemize}

	At this point, we are equipped to state the iterative method ground on Tseng's forward-backward-forward algorithm with extrapolation from the past and penalty scheme and its convergent result.
	
	\begin{tcolorbox}[colback=white]
		\textbf{Algorithm 4}:
		\begin{equation}
			\begin{aligned} 
				\mbox{Initialization}&:\; \mbox{Choose}\;  \;(x_0,v_{1,0},\dots,v_{m,0}), (y_{-1},q_{1,-1},\dots,q_{m,-1}) \in \mathcal{H} \times \mathcal{G}_1, \times \dots \times \mathcal{G}_m \\ 
				&\quad \;  \mbox{with} \;  y_{-1}=x_{0},q_{i,-1}=v_{i,0} \;\forall i = 1,\dots,m. \\
				\mbox{For}\; n\in\mathbb{N}\cup\{0\}&, \\
				\; \mbox{set}
					&:\; y_n = prox_{\lambda_n f}[ x_n-\lambda_n (\nabla h(y_{n-1}) + \sum_{i=1}^{m} L^*_i (q_{i,n-1}) )- \lambda_n\beta_n \nabla \Psi(y_{n-1})];\\
					&\quad  q_{i,n} = prox_{\lambda_n g_i^{*}} [v_{i,n} + \lambda_n ( L_i (y_{n-1}) ], \; i=1,\dots,m, \\ 
					&\quad x_{n+1}\!=\!\lambda_n\beta_n(\nabla \Psi(y_{n-1})\!-\!\nabla \Psi (y_n)) \!+\! \lambda_n (\nabla h(y_{n-1})\!-\! \nabla h(y_n))\!\\
					&\quad\quad\quad\quad+\! \lambda_n \sum_{i=1}^{m} L^*_i(q_{i,n-1}-q_{i,n})\!+\!y_n;\\
					&\quad v_{i,n+1} = \lambda_n L_i(y_n-y_{n-1}) 
					  + q_{i,n}, \; i=1,\dots,m, \\ 
			\end{aligned}
		\end{equation}
	\end{tcolorbox}
	 \noindent where  $(\lambda_n)_{n\in\mathbb{N}\cup\{0\}}$ and  $(\beta_n)_{ { n\in\mathbb{N}\cup\{0\} } }$  are positive real numbers. To analyze the convergence theorem, we need to assume the hypotheses similar to those in [\cite{2014BotCsetnek}, Section 3], which are as follows:
	 
	 $(H^{opt})\begin{cases}
	 	(i) &\partial f + N_C  \;\mbox{is maximally monotone and} \; (\ref{Convex-Opt Promblem}) \; \mbox{has an optimal solution};\\
	 	(ii) &\mbox{For every } p\in ran (N_C),  \sum_{n\in\mathbb{N}\cup\{0\} } \lambda_n\beta_n [\Psi^{*} (\frac{p}{\beta_n})-\sigma_C(\frac{p}{\beta_n})] < +\infty;\\
	 	(iii) &(\lambda_n)_{ n\in\mathbb{N}\cup\{0\} }  \in \ell^2 \setminus \ell^1.
	 \end{cases}$
	
	For the hypothesis $(H^{opt})$, some observations have already been presented in [\cite{2014BotCsetnek}], which are as follows:
		\begin{remark}\label{Remak extend hypothesis}\
			\begin{enumerate}
				\item Due to [\cite{BC-Book}, Corollary 24.4], $\partial f + N_C$ is maximally monotone, if $0\in \operatorname{sqri} (\operatorname{dom} f - C)$. This condition is fulfilled if the conditions in [\cite{BC-Book}, Proposition 6.19, 15.24] hold, for instance, $0\in \operatorname{int} (\operatorname{dom}f -C )$ or $f$ is continuous at a point in $\operatorname{dom} f \cap C$ or $\operatorname{int} C \cap \operatorname{dom} f \neq \emptyset$. 
				
				\item The condition $(ii)$ of $(H^{opt})$ is similar to the condition $(ii)$ of $(H)$ which implies to the second contion of $(H_{fitz})$ as we mentined in the Section \ref{Motivation}. Hence, we can apply our proposed convergent results in Section \ref{forward-backward-forward Penalty Schemes} for this context.

			\end{enumerate}
		\end{remark}

	With the same effort as in Section \ref{The FBF algorithm with extrapolation penalty scheme for the problem involving linearly composed and parallel-sum type monotone operators}, we may continue and provide the convergence results for \textbf{Algorithm 4} as shown below.

		\begin{theorem}\label{Thm for Convex optimization}
	Let $(x_n)_{n\in\mathbb{N}\cup\{0\}}$ be the sequences generated by \textbf{Algorithm 4} and $(z_n)_{n\in\mathbb{N}\cup\{0\}}$ the sequence defined in (\ref{weigthed_averages}). If $(H^{opt})$ and (\ref{sqri-qualification condition}) are fulfilled and  $\limsup_{n\rightarrow +\infty} \left( \frac{\lambda_n \beta_n}{\mu} + \lambda_n \beta \right) < \frac{1}{2}$, where
	$$\beta = \nu
	+ \sqrt{\sum_{i=1}^{m}  \|L_i\|^2},$$
	then $(z_n)_{n\in\mathbb{N}\cup\{0\}}$ converges weakly to an optimal solution to (\ref{Convex-Opt Promblem}) as $n\rightarrow +\infty$. If, additionally, $f$ and $g_i^{*}, i =1,\dots,m$ are strongly convex, then $(x_n)_{n\in\mathbb{N}\cup\{0\}}$ converges strongly to the unique optimal solution of the problem in (\ref{Convex-Opt Promblem}) as $n\rightarrow +\infty$.
	\end{theorem}

	\begin{remark}\
		\begin{enumerate}
			\item According to [\cite{BC-Book}, Proposition 17.10, Theorem 18.15], $g^*$ is strongly convex whenever  $g: \mathcal{H}\rightarrow \mathbb{R}$ is differentiable with Lipschitz continuous gradient.
			\item If $\Psi(x) = 0$ for all $x\in\mathcal{H}$, then \textbf{Algorithm 4} reduces to the error-free version of the iterative scheme proposed in [\cite{Tongnoi2022}] where all of the variable matrices involved are identity matrices. This scheme is used to solve the convex minimization problem given by
			\begin{align*}
				\inf_{x\in\mathcal{H}} \bigg\{ f(x) + \sum_{i=1}^{m}g_i
				(L_i x) + h(x) \bigg\}.
			\end{align*}
			
		\end{enumerate}
	\end{remark}


\section{A Numerical Experiment in TV-Based Image Inpainting}\label{A Numerical Experiment in TV-Based Image Inpainting}


In this section, we demonstrate the application  of \textbf{Algorithm 4} for solving an image inpainting problem, which involves recovering lost information in an image. The computations presented in this part were carried out with Python (version 3.7.9) on a Windows desktop computer powered by an Intel(R) Core(TM) i5-8250U processor that operated at speeds between 1.6 GHz and 1.8  GHz, and was equipped with 8.00 GB of RAM. We represent images of size $M\times N$ as vectors $x\in \mathbb{R}^{\bar{n}}$, where $\bar{n} = M\cdot N$, and each pixel denoted by $x_{i,j}, 1 \leq i \leq M, 1 \leq j \leq N$, takes values in the closed interval from $0$ (pure black) to $1$ (pure white). Given an image $b\in\mathbb{R}^{\bar{n}}$ with missing pixels (which are set to black in this case), we define $P\in \mathbb{R}^{\bar{n}\times \bar{n}}$ as the diagonal matrix with $P_{i,i}=0$ if the $i$-th pixel in the noisy image $b$ is missing, and $P_{i,i} = 1$ otherwise, for $i = 1,\dots,\bar{n}$ (noting that $\|P\| =1$). The original image is reconstructed by solving the following TV-regularized model:
	\begin{align}\label{TV model}
	\inf \{ TV_{iso}(x) : Px=b, x\in[0,1]^{\bar{n}} \}.
	\end{align}

The function $TV_{iso}:\mathbb{R}^{\bar{n}}\rightarrow \mathbb{R}$, which we use as our objective function, is defined as the isotropic total variation:
	\begin{align*}
	TV_{iso} (x) = \sum_{i=1}^{M-1} \sum_{j=1}^{N-1} \sqrt{(x_{i+1,j} - x_{i,j})^2 + (x_{i,j+1} - x_{i,j})^2} + \sum_{i=1}^{M-1} | x_{i+1,N} - x_{i,N}| + \sum_{j=1}^{N-1} |x_{M,j+1} - x_{M,j}|.
	\end{align*}

It is possible to express the problem presented in (\ref{TV model}) as an optimization problem of the form of \textbf{Problem 4} in (\ref{Convex-Opt Promblem}). This can be achieved by introducing the set $\mathcal{Y}=\mathbb{R}^{\bar{n}}\times\mathbb{R}^{\bar{n}}$, and defining the linear operator $L:\mathbb{R}^{\bar{n}}\rightarrow \mathcal{Y}$ as $L(x_{i,j}) = (L_1 x_{i,j},L_2 x_{i,j})$,where
	\begin{align*}
	L_1 x_{i,j}= 
	\begin{cases} x_{i+1,j} - x_{i,j} \quad &\mbox{if}\; i< M, \\ 
	0 &\mbox{if}\; i = M\end{cases} \quad 
	\mbox{and} 
	\quad L_2 x_{i,j} = 
	\begin{cases}
	x_{i,j+1} - x_{i,j} \quad &\mbox{if}\; j< N,\\
	0   &\mbox{if}\; j = N.
	\end{cases}
	\end{align*} 
	The operator $L$ represents a discretization of the gradient in the horizontal and vertical directions. It is worth noting that $\| L\|^2 < 8$, and its adjoint $L^*: \mathcal{Y}\rightarrow\mathbb{R}^{\bar{n}}$ is as easy to implement as $L$ itself (see [\cite{Chambolle2004,Tongnoi2022}]). Furthermore, the inner product on $\mathcal{Y}$, given by $\inner{(y,z)}{(p,q)} = \sum_{i=1}^{M} \sum_{j=1}^{N} (y_{i,j} p_{i,j} + z_{i,j} q_{i,j})$, induces a norm defined as $\|(y,z)\|_{\times} = \sum_{i=1}^{M} \sum_{j=1}^{N} \sqrt{y_{i,j}^2 + z_{i,j}^2}$ for $(y,z),(p,q) \in\mathcal{Y}$. It can be shown that $TV_{iso}(x) = \| Lx \|_{\times}$ for every $x\in\mathbb{R}^{\bar{n}}$.

	Moreover, by considering the function $\Psi:\mathbb{R}^{\bar{n}}\rightarrow\mathbb{R}, \Psi (x) = \frac{1}{2} \| Px - b\|^2$, problem  in (\ref{TV model}) can be reformulated as 
	\begin{align}\label{in practice}
	\inf_{x\in\arg\min \Psi} \{ f(x) + g_1 (Lx) \},
	\end{align}
	where $f: \mathbb{R}^{\bar{n}}\rightarrow\bar{\mathbb{R}}, f = \delta_{[0,1]^{\bar{n}}}$, $g_1: \mathcal{Y}\rightarrow \mathbb{R}, g_1(y_1,y_2) = \|(y_1,y_2) \|_{\times}$. The problem given in (\ref{in practice}) can be written as \textbf{Problem 4} (\ref{Convex-Opt Promblem}) when we set $m=1, L_1=L$, 
	 and $h=0$. It is worth noting that $\nabla \Psi(x) = P(Px-b) = P(x-b)$ for all $x\in\mathbb{R}^{\bar{n}}$, which makes $\nabla \Psi$ Lipschitz continuous with a Lipschitz constant of $\mu=1$. Therefore, the iterative scheme in \textbf{Algorithm 4} can be written as follows for every $n\geq 0$ in this particular case:
	\begin{equation}
		\left[
			\begin{aligned} 
					& y_n = prox_{\lambda_n f}[ x_n-\lambda_n L^* (q_{1,n-1}) )- \lambda_n\beta_n  P(y_{n-1}-b)],\\
					& q_{1,n} = prox_{\lambda_n g_1^{*}} [v_{1,n} + \lambda_n  L (y_{n-1}) ],\\
					& x_{n+1}\!=\!\lambda_n\beta_n  P(y_{n-1} \!-\!y_n)  +\! \lambda_n  L^* (q_{1,n-1}-q_{1,n})\!+\!y_n,\\
					&v_{1,n+1} = \lambda_n L(y_n-y_{n-1}) + q_{i,n}, 
			\end{aligned}
		\right.
	\end{equation}
		
		To implement this iterative method, we need the following formulas:
	\begin{align*}
	prox_{\gamma f}(x)= proj_{[0,1]^n} (x) \quad \forall \gamma > 0 \quad\mbox{and}\quad \forall x\in\mathbb{R}^{\bar{n}}, 
	\end{align*}
	and
	\begin{align*}
	prox_{\gamma g^{*}_1} (p,q) = proj_{S}(p,q) \quad \forall \gamma > 0 \quad\mbox{and}\quad \forall(p,q)\in \mathcal{Y},
	\end{align*}
	where $S = \{ (p,q) \in \mathcal{Y} : \displaystyle\max_{\substack{1\leq i \leq M \\ 1 \leq j \leq N}} \sqrt{p_{i,j}^2 + q_{i,j}^2} \leq 1 \}$.
	The projection operator $proj_{S} : \mathcal{Y} \rightarrow S$ is defined via 
	\begin{align*}
	(p_{i,j},q_{i,j}) \mapsto \frac{(p_{i,j}, q_{i,j})}{\max\{ 1, \sqrt{p_{i,j}^2 + q_{i,j}^2} \}},\; 1\leq i \leq M, \; 1\leq j \leq N.
	\end{align*}
	The quality of the reconstructed images was compared using the Improvement in Signal-to-Noise Ratio (ISNR), defined as follows:
	\begin{align*}
	\mbox{ISNR}(n)= 10log_{10} \left(\frac{\|x-b\|^2}{\|x-x_n\|^2} \right),	
	\end{align*}
	where $x,b$ and $x_n$ denote the original, the image with missing pixels and the recovered image at iteration $n$, respectively.
	
\begin{figure}[!htb]
		\vspace{6mm}
		\centering
		\includegraphics[width=\textwidth]
		{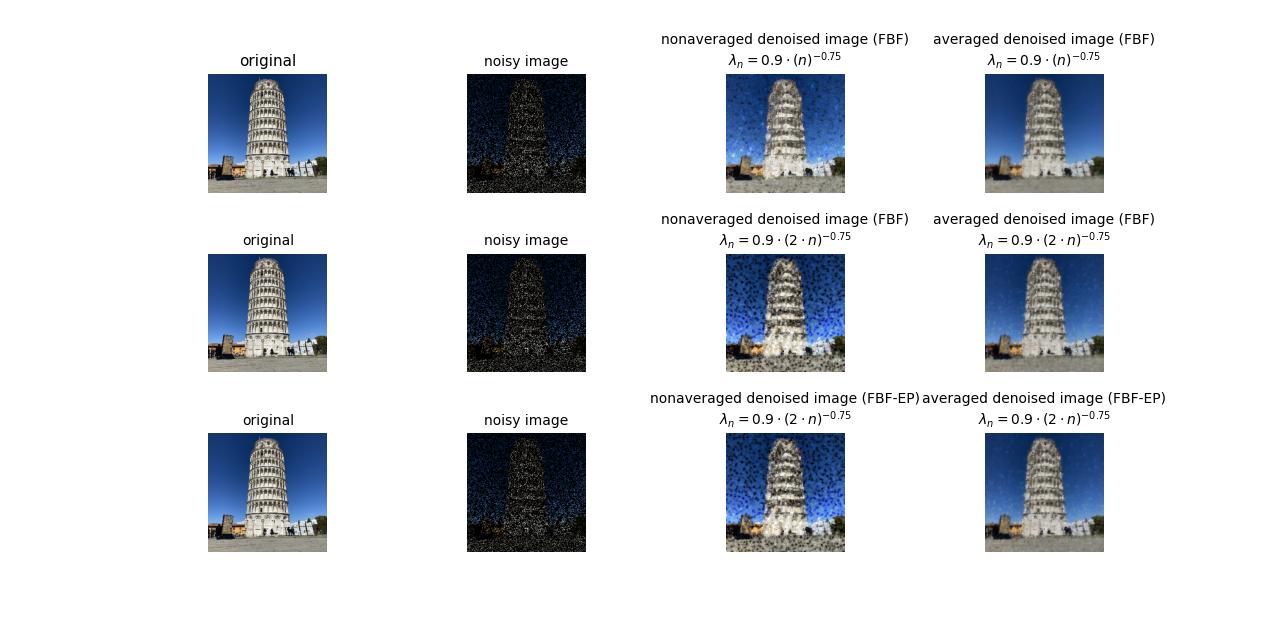}
		\caption{TV image inpainting was performed on four images, including the original image, an image with $80\%$ of its pixels missing, the non-averaged reconstructed image denoted as $x_n$, and the reconstructed image $z_n$ obtained after 2000 iterations using FBF and FBF-EP in the first and second row, respectively.}
		\label{fig:shoe pics}
		\vspace{-4mm}
	\end{figure}

	
	We evaluated the performance of the algorithm on a $256\times256$ pixel image of the Pisa Tower, using the parameters $\lambda_n = 0.9 \cdot  n^{-0.75}$, $\lambda_n = 0.9 \cdot (2\cdot n)^{-0.75}$ and $\beta_n = n^{0.75}$ for all $n \in \mathbb{N}$. The original image, as well as an image with $80\%$ randomly blacked-out pixels, were used in the test. Figure \ref{fig:shoe pics} displays the original image, the noisy image, the non-averaged reconstructed image $x_n$, and the averaged reconstructed image $z_n$ after 2000 iterations using both the Forward-Backward-Forward (FBF) in the penalty scheme with $\lambda_n = 0.9 \cdot  n^{-0.75}$, FBF in penalty scheme with $\lambda_n = 0.9 \cdot  (2\cdot n)^{-0.75}$ and Forward-Backward-Forward with Extrapolation from the Past(FBF-EP) in penalty scheme with $\lambda_n = 0.9 \cdot  (2\cdot n)^{-0.75}$. Note that $\lambda_n = 0.9 \cdot n^{-0.75}$ can not be used for FBF-EP in penalty scheme because there is no theorem to support its convergence.


\begin{table}[ht]

\setlength\extrarowheight{0pt} 
	
\begin{tabular}{@{} ccccccc @{}} 
 \toprule 
  & \multicolumn{2}{c}{FBF} & \multicolumn{2}{c@{}}{FBF} & \multicolumn{2}{c@{}}{FBF-EP}\\
  & \multicolumn{2}{c}{$\lambda_n =0.9n^{-0.75}$} & \multicolumn{2}{c@{}}{$\lambda_n =0.9(2\cdot n)^{-0.75}$} & \multicolumn{2}{c@{}}{$\lambda_n =0.9(2\cdot n)^{-0.75}$}\\
 \cmidrule(lr){2-3} \cmidrule(l){4-5} \cmidrule(l){6-7}
            & average & non-average & average & non-average & average & non-average\\
 \midrule 
 ISNR     &  11.35073   &  10.80701       & 11.32596 & 8.80449 	& 11.35116 &  8.79544\\
 CPU-time (sec) &  \multicolumn{2}{c}{94.58079} &  \multicolumn{2}{c}{95.42050}  & \multicolumn{2}{c}{80.28986}\\
 \bottomrule 
\end{tabular}
\caption{The table presents the results of ISNR and CPU time (in seconds) for both averaged and non-averaged reconstructed images obtained through 2000 iterations of the FBF and FBF-EP methods.}
\label{table:reslut of experiment}
\end{table}

\begin{figure}[!htb]
			\centering
			\includegraphics[width=13cm,height=10cm]{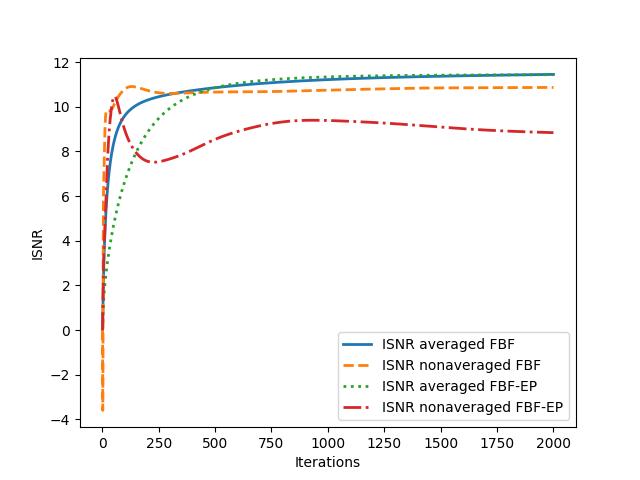}
			\caption{The figure illustrates the ISNR curves for both the averaged and non-averaged reconstructed images using the FBF  method with $\lambda_n =0.9n^{-0.75}$ and FBF-EP with $\lambda_n =0.9(2 \cdot n)^{-0.75}$ method.}
			\label{fig:graph1}
			\vspace{-2mm}
		\end{figure}

\begin{figure}[!htb]
			\centering
			\includegraphics[width=13cm,height=10cm]{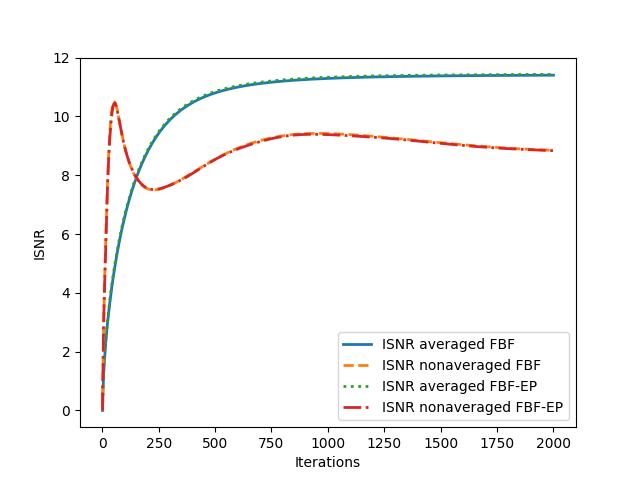}
			\caption{The figure illustrates the ISNR curves for both the averaged and non-averaged reconstructed images using the FBF  and FBF-EP methods with the same value of $\lambda_n =0.9(2 \cdot n)^{-0.75}$.}
			\label{fig:graph2}
			\vspace{-2mm}
		\end{figure}

		Figure \ref{fig:graph1} and \ref{fig:graph2} depict the evolution of the ISNR values for both the averaged and non-averaged reconstructed images using FBF and FBF-EP algorithms (in both two values of $\lambda_n$ in FBF, respectively). The theoretical outcomes concerning the sequences involved in Theorem \ref{Thm for Convex optimization} are illustrated in the figure, indicating that the averaged sequence exhibits better convergence properties than the non-averaged sequence. We notice that the behavior of ISNR values in both algorithms carry out similarly when we use the same value of $\lambda_n$ and the ISNR values of averaged FBF-EP have a bit outperform at 2000 iteration, even though we choose a bigger $\lambda_n$ in FBF for comparison. Table \ref{table:reslut of experiment} reveals that the FBF-EP approach produced the best averaged reconstructed image with an ISNR value of  11.35116, while it was 11.32596 for the averaged case of the FBF with $\lambda_n =0.9(2 \cdot n)^{-0.75}$ and 11.35073 for the averaged case of the FBF with $\lambda_n =0.9n^{-0.75}$ . However, the non-averaged reconstructed image using the FBF-EP method had the lowest ISNR value of 8.79544. Furthermore, the FBF-EP method was faster than the FBF method by around 15 seconds (compare with both two various of $\lambda_n$ in FBF), indicating that it can offer time-saving benefits in image reconstruction.\\

		\begin{remark}\
			\begin{enumerate}
			\item During the computation, the construction of weighted averaged sequences, as defined in (\ref{weigthed_averages}), requires the sequence $(x_n)_{n\in\mathbb{N}\cup\{0\}}$ in every iterative step. Thus, we computed them simultaneously in our program and measured the time consumption.
			\item The Python code is uploaded to the repository in GitHub. The reader can check it at the following link: \url{https://github.com/BurisT/inpainting_python}.
			\end{enumerate}
		\end{remark}

	\noindent\textbf{Acknowledgment.} The Royal Government of Thailand scholarship through the Development and Promotion of Science and Technology Talents Project (DPST) provided funding for this work. The author is highly appreciative of Dr. habil. Ernö Robert Csetnek's diligent direction.


\nocite{*}

\clearpage


\begin{thebibliography}{99}
	\bibitem{Attouch2010Czarnecki}
 	 Attouch, H., Czarnecki, M.-O.: Asymptotic behavior of coupled dynamical systems with multiscale
 	aspects. J. Differ. Equ. 248(6), 1315--1344 (2010)
 	
 	\bibitem{AttouchCzarnecki2011}
	 Attouch, H., Czarnecki, M.-O., Peypouquet, J.: Coupling forward-backward with penalty schemes and parallel splitting for constrained variational inequalities. SIAM J. Optim. 21, 1251--1274 (2011)
	 
 	\bibitem{Attouch2011Czarnecki} 
	 Attouch, H., Czarnecki, M.-O., Peypouquet, J.: Prox-penalization and splitting methods for constrained
variational problems. SIAM J. Optim. 21, 149--173 (2011)

	
 	
	\bibitem{BC-Book}
	Bauschke H. H.,  and Combettes P. L.:  Convex Analysis and Monotone Operator Theory in Hilbert Spaces. Springer, New York (2011)
	
	\bibitem{Bauschke2006}
	 Bauschke, H.H., McLaren, D.A., Sendov, H.S.: Fitzpatrick functions: inequalities, examples and remarks
on a problem by S. Fitzpatrick. J. Convex Anal. 13, 499--523 (2006)

	\bibitem{BanertBot2015}
	Banert, S., and Boţ, R. I.: Backward penalty schemes for monotone inclusion problems. Journal of Optimization Theory and Applications, 166(3), 930--948 (2015)

	\bibitem{BöhmSedlmayerCsetnekBot2022}
 	 Böhm A., Sedlmayer M., Csetnek E. R., and Boţ R. I.: Two Steps at a Time—Taking GAN Training in Stride with Tseng’s Method. SIAM Journal on Mathematics of
 	Data Science, 4(2), 750--771 (2022)
 	
	\bibitem{2014BotCsetnek} 
	 Boţ, R. I., and Csetnek, E. R.:  A Tseng’s type penalty scheme for solving inclusion problems involving linearly composed and parallel-sum type monotone operators. Vietnam Journal of Mathematics, 42(4), 451--465 (2014)
	 
	 \bibitem{BotCsetnek2014}
	 Boţ, R. I., and Csetnek, E. R.:  Forward-backward and Tseng’s type penalty schemes for monotone inclusion problems. Set-Valued and Variational Analysis, 22, 313--331 (2014)
	 
	\bibitem{BotCsetnekNimana2018}
 	Boţ, R. I., Csetnek, E. R., and Nimana, N.: Gradient-type penalty method with inertial effects for solving constrained convex optimization problems with smooth data. Optimization Letters, 12, 17--33 (2018)
 	
 	 
	
	 \bibitem{Chambolle2004}
	 Chambolle, A.: An algorithm for total variation minimization and applications. J. Math. Imaging Vis. 20, 89--97 (2004)
	
	\bibitem{CombettesPesquet2012}
	 Combettes, P.L., Pesquet, J.-C.: Primal-dual splitting algorithm for solving inclusions with mixtures of composite, Lipschitzian, and parallel-sum type monotone operators. Set-Valued Var. Anal. 20, 307--330 (2012)
	 
	\bibitem{Fitzpatrick1988}
	 Fitzpatrick, S.: Representing monotone operators by convex functions. In: Workshop/Miniconference
	on Functional Analysis and Optimization, Canberra 1988. Proceedings of the Centre for Mathematical
	Analysis, vol. 20, pp. 59--65. Australian National University, Canberra (1988)

	 \bibitem{Noun2013Peypouquet} 
 	Noun, N., Peypouquet, J.: Forward–backward penalty scheme for constrained convex minimization with-out inf-compactness. J. Optim. Theory Appl. 158, 787--795 (2013)
 	
 	\bibitem{Malitsky2020Tam}
 	 Malitsky Y., and Tam M.K.: A forward-backward splitting method for monotone
 	inclusions without cocoercivity. SIAM Journal on Optimization, 30(2), 1451--1472 (2020)
 	
 	\bibitem{Peypouquet2012}
 	Peypouquet, J.: Coupling the gradient method with a general exterior penalization scheme for convex
 	minimization. J. Optim. Theory Appl. 153, 123--138 (2012)
	 
 	\bibitem{Popov1980} Popov L. D.: A modification of the Arrow-Hurwicz method for search of saddle points. Mathematical notes of the Academy of Sciences of the USSR, 28(5), 845--848 (1980)
 		
 	\bibitem{Rakhlin2013Sridharan}
 	 Rakhlin A., and Sridharan K.: Online learning with predictable sequences. In Conference on Learning Theory (pp. 993--1019). PMLR. (2013, June)
 	
 	\bibitem{RakhlinSridharan2013}
 	Rakhlin S., and Sridharan K.: Optimization, learning, and games with predictable sequences. Advances in Neural Information Processing Systems, 26 (2013)
 	
	\bibitem{rockafellar1997convex}
	Rockafellar, R. T.: Convex Analysis. Princeton University Press. (1970)

	
	
	\bibitem{Rockafellar1970}
	 Rockafellar, R.T.: On the maximal monotonicity of subdifferential mappings. Pac. J. Math. 33, 209--216 (1970)
	

 
	\bibitem{Tongnoi2022}
	Tongnoi, B.: Tseng’s Algorithm with Extrapolation from the past Endowed with Variable Metrics and Error Terms. Numerical Functional Analysis and Optimization, 1--37 (2022).
	
 	
 	\bibitem{Tseng2000}
 	 Tseng, P.: A modified forward-backward splitting method for maximal monotone mappings. SIAM J. Control. Optim. 38(2), 431--446 (2000)
 	
 	
 	
 	
 	\bibitem{Yang2022LiLan}
	 Yang, S., Li, X., and Lan, G.: Data-Driven Minimax Optimization with Expectation Constraints. arXiv preprint arXiv:2202.07868 (2022)
	
 	
 



	
\end{thebibliography}

\afterpage{
}
\end{document}